%% file: NM_homomorphisms_of_Beurling_Algebras_3.tex
\documentclass[11pt]{article}

\usepackage{amsmath,amsfonts,amssymb,amsthm,mathtools,bbm, geometry}

\input{Macros}

\usepackage[all]{xy}

\def\S{{\cal S}} 
\def\omg{{\om_G}} 
\def\omf{{\om_F}}
\def\omhk{{\om_{H/K}}}
\def\omx{{\om_\times}}
\def\omgi{{\om_G^{-1}}} 
\def\omfi{{\om_F^{-1}}}
\def\omi{{\om^{-1}}}
\def\omhki{{\om_{H/K}^{-1}}}
\def\C{\mathbb{C}}
\def\R{\mathbb{R}}
\def\s{{\rm s}}
\def\T{\mathbb{T}} 
\def\Cx{\C^\times}
\def\Jgt{J_{\gamma, \theta}} 
\def\jgt{j_{\gamma, \theta}}

\begin{document}

\title{Norm-Multiplicative Homomorphisms of Beurling  Algebras  }
\author{Matthew E. Kroeker\footnote{This author was  supported by a 2018 NSERC USRA.}, Alexander Stephens$^*$, Ross Stokke\footnote{
This author was partially supported   by an NSERC grant. }  \   and  Randy Yee\footnote{This author was supported by R. Stokke's NSERC grant in the summer of 2013.}   }
\date{}
\maketitle

\begin{abstract}{\small   We introduce and study ``norm-multiplicative" homomorphisms $\vp: \L^1(F) \ra \M_r(G)$   between group and measure algebras,  and $\vp: \L^1(\omf) \ra \M(\omg)$ between Beurling group and measure algebras, where $F$ and $G$  are locally compact groups with continuous weights  $\omf$ and $\omg$.  Through a unified approach we recover, and sometimes strengthen, many of the main known results concerning homomorphisms and isomorphisms between these (Beurling) group and measure algebras.  We provide a first description of all positive homomorphisms $\vp: \L^1(F) \ra \M_r(G)$.  We state  versions of  our results that describe a variety of  (possibly unbounded) homomorphisms $\vp: \C F \ra \C G$ for (discrete) groups $F$ and $G$. 

\smallskip

\noindent{\em Primary MSC codes:}  43A20, 43A10, 43A15, 43A22,  16S34   

 }
\end{abstract} 

\section{Introduction} 

Throughout this article, $(F, \omf)$ and $(G, \omg)$ are weighted locally compact groups. The Beurling group and measure algebras over $(F, \omf)$ are denoted by $\L^1(\omf)$ and $\M(\omf)$ respectively; when $\omf$ is the trivial weight,  we obtain the usual group and measure algebras, denoted herein by $\L^1(F)$ and $\M_r(F)$.  

A natural endeavour in any area of mathematics is to describe the structure preserving maps between like objects. The existence of an algebra isomorphism between $\C\, \mathbb{Z}_4$ and $\C (\mathbb{Z}_2 \times \mathbb{Z}_2)$ shows that the algebraic structures of group and measure algebras do not alone determine the underlying group \cite{Kal-Woo}.  In abstract harmonic analysis, a well-studied question thus asks for a description  of isomorphisms   between  various Banach algebras built over  locally compact groups that take account of the additional structure of these algebras, and whether the existence of such isomorphisms determine the underlying locally compact groups.  Work on this problem began with the study of isometric and bipositive isomorphisms between group algebras and measure algebras  by Kawada, Wendel,  Johnson and Strichartz in  \cite{Kaw, Wen1, Wen2, Joh, Str} (and elsewhere). These authors showed that the Banach algebra and order structures of group and measure algebras determine their underlying locally compact groups, i.e. that $\L^1(F)$ and $\M_r(F)$ are complete invariants for the locally compact group $F$. The literature concerned with determining whether other types of Banach algebras over $F$ are complete invariants for $F$ is quite extensive; a small sample of relevant papers is \cite{Gha-Lau1,Gha-Lau-Lo, Gha-Lau2, Lau-McK, Wal}.  

More generally, the ``homomorphism problem" in abstract harmonic analysis asks for a description of all homomorphisms $\vp: \L^1(F) \ra \M_r(G)$, (the dual problem asks for a description of all homomorphisms between Fourier and Fourier--Stieltjes algebras). This problem was completely solved in the abelian case  by Paul Cohen in \cite{Coh}, but remains open in general. Contractive homomorphisms $\vp: \L^1(F) \ra \M_r(G)$ were first described by Greenleaf in \cite{Gre}, and a characterization of such homomorphisms via a Cohen-type factorization is found in \cite{Sto1}; see \cite{Sto1} for more about the history of this old problem. There  had been no known description of the positive (order-preserving) homomorphisms $\vp: \L^1(F) \ra \M_r(G)$.

Beurling group and measure algebras provide the setting for a rather active area of research, e.g., see    \cite{Dal-Lau, Fil-Sal, Gha1,Gha2, Kan, Neu, Rei-Ste, Sam, She-Zha}.    
 Zadeh \cite{Zad} and Ghahramani-Zadeh \cite{Gha-Zad} recently showed that the Beurling group algebra $\L^1( \omf)$   and the  Beurling measure algebra $\M(\omf)$ are,  up to an isometric/bipositive isomorphism,   complete invariants for the underlying weighted locally compact group $(F, \omf)$.

Extending  work of Pym \cite{Pym} and Greenleaf \cite{Gre} to the setting of Beurling measure algebras, in Section 2 we describe all norm-one, and positive norm-one, idempotents in $\M(\omg)$. In Section 3 we introduce and completely describe the norm-multiplicative subgroups --- or  $\NM_\omg$-subgroups ($\NM$-subgroups in the case of the trivial weight) ---   of $\M(\omg)$: these are the  (convolution) subgroups $\Gamma$ of $\M(\omg)$ on which the norm $\| \cdot \|_\omg: \Gamma \ra (0, \infty)$ is a homomorphism. The collection of $\NM_\omg$-subgroups of $\M(\omg)$ properly contains the collection of all contractive subgroups of $\M(\omg)$, so the results in this section extend Greenleaf's description from \cite{Gre} (and \cite{Sto1}, which contains  the topological description)  of the contractive subgroups of $\M_r(G)$. (Note that we rely heavily on the results and methods from \cite{Gre, Sto1} in Section 3.)  Every positive subgroup of $\M(\omg)$ with identity in $\M_r(G)$ is also a $\NM$-subgroup of  
$\M_r(G)$, and we thus obtain a first description of all positive subgroups of $\M_r(G)$. 

We call a bounded homomorphism $\vp : \L^1(\omf) \ra \M(\omg)$ a $\NM_\omg$- ($\NM$-)homomorphism if the image of $F$ under the canonical strict-to-weak$^*$ continuous extension of $\vp$ to $\M(\omf)$  is a $\NM_\omg$-subgroup of $\M(\omg)$ ($\NM$-subgroup of $\M_r(G)$). We introduce some ``basic homomorphisms" and characterize all $\NM_\omg$-homomorphisms (and $\NM$-homomorphisms for which $\|\vp(\delta_{e_G})\|_\omg =1$)  as those that factor into a particular product of three basic homomorphisms.  While $\NM_\omg$- and $\NM$-homomorphisms $\vp: \L^1(\omf) \ra \M(\omg)$ may not be contractive, we observe that the $\NM$-homomorphisms $\vp: \L^1(F) \ra \M_r(G)$ are precisely the contractive homomorphisms and recover Theorem 5.11 of \cite{Sto1} (though we make significant use herein of \cite{Sto1}).  We describe all positive homomorphisms $\vp: \L^1(\omf) \ra \M(\omg)$ such that $\|\vp(\delta_{e_G})\|_\omg =1$, and thereby describe \it all \rm positive homomorphisms $\vp: \L^1(F) \ra \M_r(G)$, also a new result. We conclude Section 4 by recovering the main results concerning isometric and bipositive isomorphisms $\vp: \L^1(\omf) \ra \L^1(\omg)$ and $\vp: \M(\omf) \ra \M(\omg)$ from \cite{Zad} and \cite{Gha-Zad}; moreover, we show that positive isomorphisms are always bipositive, thereby strengthening the results from \cite{Kaw, Gha-Zad}. This paper provides a unified approach to obtaining many of the main known results about isomorphisms and homomorphisms between (Beurling) group and measure algebras. 

In Section 5, we state  versions of some our  main results describing various (not necessarily bounded) homomorphisms $\vp: \C F \ra \C G$. We believe these results are also new and are hopeful that they will  be of interest to algebraists working in group rings.

\subsection*{Preliminary results and notation}

Throughout this article, $G$ is a locally compact group and $\om: G \ra (0,\infty)$ is a continuous \it weight function \rm  satisfying
$$\om(st) \leq \om(s)\om(t) \ \ (s, t \in G) \qquad {\rm and} \qquad \om(e_G) = 1;$$
the pair $(G, \om)$ is called a \it weighted locally compact group. \rm    When dealing with two weighted locally compact groups $F$ and $G$, we use $\omf$ and $\omg$ to specify their weights.   Let $\lambda$ denote a fixed left Haar measure on $G$, with respect to which the group algebra $\L^1(G)$ and $\L^\infty(G)= \L^1(G)^*$ are defined in the usual way.  The Beurling group algebra   $\L^1(G, \om)$,  or simply $\L^1(\om)$,  is composed of all functions $f$ such that   $ \om f$ belongs to $\L^1(G)$, with $\|f\|_{1,\om} := \|\om f\|_1$ and convolution product. If   ${\cal S}(G)$   is a closed subspace of $\L^\infty(G)$,  $\psi \in {\cal S}( \om^{-1})$ exactly when ${\psi / \om} \in {\cal S}(G)$;  putting $\|\psi \|_{\infty, \om^{-1}} = \| \psi / \om\|_\infty$,  ${\cal S}( \om^{-1}) $ is a Banach space and 
$ S : {\cal S}( \om^{-1})  \ra {\cal S}(G): \psi \mapsto \psi /  \om$ is an isometric linear isomorphism. 

Following \cite{Fel-Dor} and \cite{Sto2}, $\fS(G)$ is the $\delta$-ring of compacted-Borel subsets of $G$, i.e., the set of all Borel subsets of $G$ with compact closure, and  $\M(G)$ denotes the space of all complex regular compacted-Borel measures $G$, i.e., all complex regular measures $\mu$ defined on $ \fS(G)$.    A measure  $\nu$ in $\M(G)$ belongs to $\M( \om)$ if  $\om \nu \in \M_r(G)$, where $\M_r(G)$ is the Banach algebra of bounded measures in $\M(G)$ with its usual norm $\| \c \|$. With the  norm $\|\nu \|_\om : = \|\om \nu \|$,  $\M(\om)$ can be identified with the dual of $C_0(\om^{-1})$ through the pairing 
$$\l \nu, \psi\r_\om = \int \psi \, d\nu \qquad (\nu \in \M(\om), \ \psi \in C_0(\om^{-1})),$$
and $\nu \mapsto \om \nu: \M(\om) \ra \M_r(G)$ is the dual map of $\phi \mapsto  \phi \om:  C_0( G) \ra C_0(\om^{-1})$; thus 
\beq \label{DualityIntegralEqn}
\int \phi \, d(\om \nu) = \int  \phi \om \, d\nu \quad \text{for}  \quad \phi \in C_0(G), \ \nu \in \M(\om).
\eeq
With respect to the convolution product, $\M(\om)$ is a dual Banach algebra, called the Beurling measure algebra. 
  For details see \cite{Fel-Dor} and  \cite{Sto2}. (The paper \cite{Sto2} corrects an earlier definition of the Beurling measure algebra found in the literature and was written specifically to validate the work herein, without requiring that we assume our weights are bounded away from zero.)

\br \rm \label{Remark mu_e}   1. Let $M(G)$ denote the usual measure algebra of complex regular measures defined on $\mathfrak{B}(G)$, the $\sigma$-algebra of Borel subsets of $G$.  For $\mu \in \M_r(G)$, $\mu_e \in M(G)$ where for $A \in \mathfrak{B}(G)$, $\mu_e(A):= \lim \mu(C)$, with the limit taken over the directed family of compact subsets $C$ of $A$.  The map $\mu \mapsto \mu_e: \M_r(G) \ra M(G)$ is an order-preserving  weak$^*$-homeomorphic isometric algebra isomorphism such that  $\int  \phi \, d \mu  = \int \phi \, d \mu_e$ for $\phi \in C_0(G)$, $|\mu_e| = |\mu|_e$, and  the supports of $\mu$ and $\mu_e$ are equal, i.e., ${\rm s}(\mu) = {\rm s}(\mu_e)$ \cite[Remarks 1.1, 1.2, 2.1]{Sto2} and \cite[II.8.15]{Fel-Dor}. As a dual Banach algebra, $\M_r(G)= C_0(G)^*$ can thus be used interchangeably with $M(G)= C_0(G)^*$.  For the sake of consistency, and because $\M_r(G)$ is more natural in this setting, we will mostly use $\M_r(G)$ in place of   $M(G)$. 

\smallskip 

\noindent 2.  When $H$ is a closed subgroup of $G$, we can  identify $\M_r(H)$ with the closed subalgebra $\M_{r,  H}(G) = \{ \mu \in \M_r(G): {\rm s}(\mu)\subseteq H \}$ of $\M_r(G)$,  its image under the weak$^*$-continuous isometric algebra isomorphism $R_H^*$, where $R_H: C_0(G) \ra C_0(H)$ is the restriction operator, a quotient map \cite[Proposition 5.3]{Sto1}. Letting $\om_H= \om\vert_H$,   $R_H$ also defines a quotient map of $C_0(\om^{-1})$ onto $C_0(\om_H^{-1})$, and $R_H^*$ is a weak$^*$-continuous isometric algebra isomorphism of $\M(\om_H)$ onto the weak$^*$-closed subalgebra $\M_H(\om) = \{\mu \in \M(\om): \s(\mu) \subseteq H\}$ of $\M(\om)$; we will often identify $\M(\om_H)$ with $\M_H(\om) \subseteq \M(\om)$ through $R_H^*$. For $\eta$ in $\M_r(H)$ or $\M(\om_H)$, $R_H^*(\eta)(A) = \eta(A \cap H)$ for $A \in \fS(G)$; if $\mu$ is  in $\M_r(G)$ or $\M(\om)$ with $\s(\mu) \subseteq H$, $\mu = R_H^*(\mu_H)$ where $\mu_H= \mu \large\vert_{\fS(H)}$. 

\smallskip 

\noindent 3.  If $\mu \in \M(G)$ and  $f \in \L^1(\mu)$, then $f$ must vanish off a $\sigma$-compact subset of $G$, a subtlety that requires some care \cite[Remark 1.2]{Sto2}. However, if $f \in \L^1(\mu)$, then $f \in \L^1(\mu_e)$ and $\int f \, d \mu = \int f \, d \mu_e$ \cite[II.8.15]{Fel-Dor}; as noted above, $\mu_e \in M(G)$ when $\mu \in \M_r(G)$. This will be used on several occasions.   

\smallskip 

\noindent 4.  Letting $\M_{cr}(G)$ denote the compactly supported measures in $\M(G)$, $\M_{cr}(G)$ is a dense subalgebra of both $\M_r(G)$ and $\M(\om)$ \cite[Remark 2.1.3]{Sto2}. 

\er


Let $\S(\om^{-1})$ be a left introverted subspace of $\L^\infty(\om^{-1}) = \L^1(\om)^*$ such that 
\beq \label{Introverted space containment assumptions}   C_0(\om^{-1}) \subseteq \S(\om^{-1}) \subseteq LUC(\om^{-1}).  \eeq
Equivalently, under the assumption (\ref{Introverted space containment assumptions}), $\S(\om^{-1})$ is a left introverted subspace of $\ell^\infty(\om^{-1}) = \ell^1(\om)^*$;  $\S(\om^{-1})^*$ is a Banach algebra with respect to its left Arens product  given by 
$$\l m \sq n, \psi \r = \l m, n \sq \psi \r \quad {\rm where} \quad  (n \sq \psi)(s) =  n (\psi \c s), \ \ (\psi \c s) (t) = \psi (st), $$
for $m, n \in \S(\om^{-1})^*, \psi \in  \S(\om^{-1}), s, t \in G$ \cite[Proposition 2.9]{Sto2}. The (left and right) Arens product and convolution agree on $\M(\om)$.  
Let $so_l^\om$ and $so_r^\om$ denote the left- and right-strict topologies on $\M(\om)$ taken with respect to the ideal $\L^1(\om)$, i.e.,   the locally convex topologies  respectively generated by the semi-norms
$p_g(\nu) = \|g * \nu \|$ and $q_g(\nu)= \| \nu * g \|$ for $g \in \L^1(\om), \nu \in
\M(\om)$,  and let $so^\om = so_l^\om \vee so_r^\om$ denote the strict topology on $\M(\om)$; when $\om \equiv 1$ is the trivial weight, we just use $so_l$, $so_r$, $so$.  As $\L^1(\om)$ has a contractive approximate identity,  (the unit ball of) $\L^1(\om)$ is $so_l/so_r$-dense in (the unit ball of) $\M(G, \om)$.   For brevity and clarity, we will sometimes denote the weak$^*$-topology on $\S(\om^{-1})^*$, including $\M(\om)$,  by $\sigma^\om$  and $\sigma$ when   $\om \equiv 1$.  Letting $\Theta: \M(\om) \hookrightarrow \S(\om^{-1})^*$ be the $so_l^\om-\sigma^\om$ (and $so_l^\om-\sigma^\om$) continuous homomorphic isometric embedding defined in \cite{Sto2}, $\Theta$ maps into the topological centre $Z_t(\S(\om^{-1})^*)$ of $\S(\om^{-1})^*$,  and extends both $\eta_\S: \L^1(\om) \hookrightarrow \S(\om^{-1})^*$ and $\eta_\S: \ell^1(\om) \hookrightarrow \S(\om^{-1})^*$ \cite[Proposition 2.9]{Sto2}. Thus, $\sigma^\om$ on $\M(\om)$ is contained in $so_l^\om$ and $so_r^\om$. When convenient, we will identify $\M(\om)$ with its copy $\Theta(\M(\om))$ in $\S(\om^{-1})^*$.

For $x \in G$, $\delta_x \in \M(\om)$ is the Dirac measure at $x$ and $\Delta_G = \{ \delta_x: x \in G\}$. We will often identify $\C G$ with the linear span of $\Delta_G$ in $\M(\om)$ or $\S(\om^{-1})^*$. For $f \in \L^1(\om)$, $f \mapsto \delta_x * f, f * \delta_x: G \ra \L^1(\om)$ are continuous (e.g., see \cite[Lemma 3.1.5]{Zad0}) and $C_{00}(G) \subseteq \S(\om^{-1})$, so it is easy to see that $x \mapsto \delta_x$ is a topological group isomorphism of $G$ onto $\Delta_G$ when $\Delta_G$ has either its relative $so_l^\om/so_r^\om $-topology from $\M(\om)$ or its $\sigma^\om$-topology from $\S(\om^{-1})^*$; thus, these topologies agree on $\Delta_G$. 

Since the $\| \c\|_1$-unit ball of $\C G$ is dense in the unit ball of $\ell^1(G)$ and $S: \ell^1(G) \ra \ell^1(\om): f \mapsto f \om^{-1}$ is a linear isometric isomorphism mapping $\C G$ onto itself, the $\| \c \|_\om$-unit ball of $\C G $ is also dense in the unit ball of  $\ell^1(\om)$.  As $\eta_\S(\ell^1(\om))_{\| \c \|_\om \leq 1}$ is $\sigma^\om$-dense in $(\S(\om^{-1})^*)_{\| \c \|\leq 1}$,  the $\| \c \|_\om$-unit ball of $\C G $  is $\sigma^\om$-dense in the unit ball of $\S(\om^{-1})^*$. Also,  the $\| \c \|_\om$-unit ball of $\C G $  is $so_l^\om$-dense  in  the unit ball of $\M(\om)$ by \cite[Lemma 3.1.9]{Zad0}, (the proof of which holds verbatim using $\M(\om)$ as defined here and in \cite{Sto2}, and either $so_r^\om$ --- used in \cite{Zad0} --- or $so_l^\om$).  These observations and the proof of \cite[Lemma 1.1]{Sto2} yield the following statement.  

\blem \label{WK* continuous Homom Lemma} Let $S(\omf^{-1})$  and $S(\omg^{-1})$ be left introverted spaces with $$C_0(\omf^{-1}) \subseteq \S(\omf^{-1}) \subseteq LUC(\omf^{-1}) \quad {\rm   and}  \quad C_0(\omg^{-1}) \subseteq \S(\omg^{-1}) \subseteq LUC(\omg^{-1}).$$ Then a $\sigma^\omf-\sigma^\omg$ continuous (resp. $so_l^\om / so_r^\om - \sigma^\omg$ continuous) linear map 
$\vp: \S(\omf^{-1})^* \ra \S(\omg^{-1})^*$ (resp. $\vp: \M(\omf) \ra \S(\omg^{-1})^*$) is uniquely determined by its values on $\Delta_F$ and if  $\vp(\Delta_F)$ is contained in  $Z_t(\S(\omg^{-1})^*)$, the topological centre of $\S(\omg^{-1})^*$, then $\vp$ is a homomorphism if and only if for each $s, t\in F$, $\vp(\delta_s * \delta_t) = \vp(\delta_s) \sq \vp (\delta_t)$. 
\elem 
 Unless explicitly stated otherwise, and in Section 5,   mappings between normed spaces are  always assumed to be bounded and  linear. Any undefined notation is found in \cite{Sto2}.

\section{Norm-one idempotents in $\M(\om)$}

Throughout this section, $(G, \om)$ is a weighted locally compact group. In the case that $\om$ is the trivial weight, the following lemma is \cite[Theorem 4.1]{Pym}. Our proof follows Pym's approach --- see section 3 and the proof of Theorem 4.1 of  \cite{Pym} --- with some adjustments to take account of the weight.  When $K$ is a compact subgroup of $G$, $m_K$ denotes the normalized  Haar measure on $K$, usually viewed as a measure in $\M_{cr}(G)$, as described in Remark \ref{Remark mu_e}.2.

\blem \label{PymLemma}  Let $\nu \in \M(\om)$ be a positive norm-one idempotent. Then $K = \s(\nu)$  is a compact subgroup of $G$ such that $\om \nu = m_K$. 
\elem 

\begin{proof}  Let $\mu = \om \nu$, so $\mu \in \M_r(G)$ and $\|\mu\| = \| \nu\|_\om=1$. Let $K = \s(\mu ) = \s(\nu)$. Take $\phi \in C_0(G)^+$. Then $\mu \c \phi \in C_0(G)^+$, where 
$$\mu \c \phi(s) = \l \mu, \phi \c s \r = \int \phi \c s \, d \mu \qquad (s \in G),$$
(e.g., by the non-weighted case of \cite[Corollary 2.6]{Sto2}). Taking $a = a_\phi\in G$ such that $$\ds \mu \c \phi(a) = \max_{x \in G} \mu \c \phi(x),$$ we \it claim \rm that  $\mu \c \phi(ax) =  \mu \c \phi(a)$ for each $x \in K$:   

Employing (\ref{DualityIntegralEqn}), Remark \ref{Remark mu_e},  \cite[Corollary 2.6 ]{Sto2} and \cite[Identity (14)]{Sto2}, we have
\beqs  \mu \c \phi(a) & = &  \int \phi \c a \, d (\om \nu) = \l \nu, (\phi \c a) \om \r_\om = \l \nu * \nu, (\phi \c a) \om \r_\om 
 =  \int \int (\phi \c a)(st) \om(st) \, d\nu(s) \, d\nu(t) \\ &  \leq &  \int \int (\phi \c a)(st) \om(s) \om( t) \, d\nu(s) \, d\nu(t)  =   \int \int t\c (\phi \c a)(s)  \, d(\om\nu)(s) \, \om(t)  d\nu(t)\\
& = &   \int  (\phi \c a)\c \mu (t)  \, d(\om\nu)(t) = \l \mu, (\phi \c a) \c \mu \r =  \l \mu, \mu \c (\phi \c a)\r\\
& = & \int \mu \c \phi(as) \, d \mu(s) = \int \mu \c \phi(as) \, d \mu_e(s)  \\
& \leq & \int \mu \c \phi(a) \, d \mu_e(s) = \mu \c \phi(a) \| \mu \| = \mu \c \phi(a).
\eeqs 
The function $s \mapsto h(s) = \mu \c \phi(a) -\mu \c \phi(as)$ is continuous, non-negative on $G$ and the above calculation shows that $\ds \int h \, d\mu_e = 0$, so $h$ vanishes on $\s(\mu_e) = \s (\mu) = K$, giving the claim.

 Taking $\phi \in C_{00}(G)^+$ such that $\l \mu, \phi \r = \delta >0$ and $a = a_\phi$ as above, 
$$\mu \c \phi(as) = \mu \c \phi(a) \geq \mu \c \phi (e_G) = \delta>0 \quad \text{for all } s \in K.$$ 
Since $\mu \c \phi \in C_0(G)$, it takes values less than $\delta$ off a compact set, from which it follows that $K$ is compact. By \cite[Proposition III.10.5]{Fel-Dor}, $K = \s(\nu) = \s(\nu * \nu) = \ov{\s(\nu) \s(\nu)} =  K^2$. Thus, $K$ is a compact cancellative semigroup, and is therefore a compact group by \cite[Theorem 1]{Num}.

Observe that $\nu \in \M_r(G) \cap \M(\om)$, since it has compact support. Letting $\nu_K = \nu\large{\vert}_{\fS(K)}$,  $\mu_K = \mu\large{\vert}_{\fS(K)}$, and $\om_K = \om\large{\vert}_K$, one can readily check that $\om_K \nu_K = (\om \nu)_K = \mu_K$, so by Remark \ref{Remark mu_e}.2,  $\|\nu_K\|_{\om_K} = \|\mu_K\| =1$ and $\nu_K * \nu_K = \nu_K$. The above claim applied to $\nu_K$ and $\mu_K$ shows that for any $\phi \in C(K)^+$, 
$\mu_K \c \phi(a) =  \mu_K \c \phi(as)$ for each $s \in K$, where $\mu_K \cdot \phi$ takes its maximum value at $a$; hence, $\mu_K \cdot \phi$ is constant on $K$. Thus, 
 $\mu_K$ is a translation-invariant positive norm-one measure on $K$, i.e., $\mu_K$ is the normalized Haar measure on $K$ (and $\mu = R_K^*(\mu_K)$ is normalized Haar measure on $K$, viewed as a measure on $G$). 
\end{proof}

\bt  \label{Positive Idempotent Thm} Let  $\nu \in \M(\om)$ be a positive norm-one idempotent. Then  there is a compact subgroup $K$ of $G$ such that $\nu = m_K$ and $\om \equiv 1$ on $K$. 
\et 

\begin{proof}  By Lemma \ref{PymLemma}, $K = \s(\nu)$ is a compact subgroup of $G$ and $\om \nu = m_K$ or, equivalently, $\nu = (1/\om) m_K$.  We can assume that $G=K$; otherwise, as shown in the proof of Lemma \ref{PymLemma}, we can work with $\nu_K = \nu \large{\vert}_{\fS(K)}$ and $\om_K = \om\large{\vert}_K$. As $m_K$ is also an idempotent in $\M_r(K) = \M_{cr}(K) = \M(\om)$, for any $\phi \in C(K)^+$, 
\beqs \l m_K, \phi/\om\r & = &  \int \phi  \, {1 \over \om} \, d m_K  =   \int \phi \, d\left({1 \over \om}  m_K\right) = \l \nu, \phi\r_\om = \l \nu * \nu, \phi \r_\om 
 \\ & =  &   \iint \phi(st)\, d \nu(s)\, d \nu(t) 
  =  \iint \phi(st)  {1\over \om(s)} \, dm_K(s) \,  {1\over \om(t)} \, dm_K(t) \\
  &  \geq & \iint {\phi\over \om}\left(st\right) \, dm_K(s) \, dm_K(t)  = \l m_K * m_K, \phi/ \om \r = \l m_K, \phi/\om \r. 
\eeqs 
Taking $\phi = \om$ (on $K$), we obtain 
$$\iint {\om(st) \over \om(s)\om(t)} \, dm_K(s) \, d m_K(t) = \int 1 \, dm_K(s) \, dm_K(t).$$
Letting $\ds g(s,t) = 1 - {\om(st) \over \om(s)\om(t)}$,   $g$ is continuous and non-negative on $K \times K$, and $$\ds \int_{K \times K} g(s,t) \, d(m_K \times m_K)(s,t) = 0.$$  Hence, $g$ vanishes on 
$K \times K$, meaning that $\om(st) = \om(s) \om(t)$ for $s,t \in K$. Thus, $\om(K)$ is a compact subgroup of $(0, \infty)$, and therefore $\om \equiv 1$ on $K$. 
\end{proof}

\br  \label{Constant on Cosets Remark Etc}   \rm  If $\om \equiv 1$ on a subgroup $K$ of $G$, then $\om$ is constant on the cosets of $K$. Indeed,   for $s \in G$ and $k\in K$,  
$\om(sk) \leq \om(s) \om(k) = \om (s) \leq  \om (sk)\om(k^{-1}) = \om (sk)$, so $\om (sk) = \om(s)$. In this case,  $\om_{G/K}(sK) := \om(s)$ is hence  a  well-defined continuous function on $G$; moreover,  $\om_{G/K}$ is a weight function on $G/K$ when $K\lhd G$.  
 \er
 
 In the case of the trivial weight, the following result contains, and gives independent proofs of,  Theorems 2.1.1 and 2.1.2 of \cite{Gre}.  In the proof, we use the fact that  $\eta \leq \sigma$ if and only if $\om \eta \leq \om \sigma$ for $\eta, \sigma \in \M(G)^+$, which is readily checked. Also, it is clear from the definition of $\mu_e$ (see Remark \ref{Remark mu_e}.1) that  $\eta \leq \sigma$ if and only if $\eta_e \leq \sigma_e$.   
 
 \bp \label{Multiplicative Tot Var Prop}   Suppose that $\mu, \nu \in \M(\om)$ and $\| \mu * \nu \|_\om = \|\mu\|_\om \|\nu\|_\om$. Then $|\mu * \nu| = |\mu | *|\nu|$ and $\s(\mu * \nu) = \overline{\s(\mu) \s(\nu)}$. 
 \ep 
 
 \begin{proof} By \cite[III.10.3]{Fel-Dor}, $|\mu * \nu| \leq |\mu | *|\nu|$. To establish the reverse inequality, we show that $\om |\mu * \nu| \leq \om |\mu | *|\nu|$. To this end, suppose towards a contradiction that $\om |\mu | * |\nu| (A) >  \om |\mu  *\nu| (A)$ for some $A \in \fS(G)$. Then 
 \beqs  \|\mu* \nu\|_\om &=& \| (\om |\mu*\nu|)_e\| = (\om |\mu*\nu|)_e(G) = (\om |\mu*\nu|)_e(A) + (\om |\mu*\nu|)_e(G\bs A)\\ & \leq &  \om |\mu*\nu| (A) + (\om |\mu|*|\nu|)_e(G\bs A) \\
 & < & \om |\mu|*| \nu| (A) + (\om |\mu|*|\nu|)_e(G\bs A) = (\om |\mu|*|\nu|)_e(G) = \| |\mu| * |\nu|\|_\om \\ 
 & \leq &  \|\mu\|_\om\|\nu\|_\om = \|\mu*\nu\|_\om. 
 \eeqs
 By \cite[Proposition III.10.5]{Fel-Dor}, we now obtain $$\s(\mu * \nu) = \s(|\mu * \nu|) = \s(|\mu| * |\nu|) = \overline{\s(|\mu|) \s(|\nu|)} = \overline{\s(\mu) \s(\nu)}. \qedhere $$
 \end{proof} 
 
 \bc \label{Multiplicative Norm Corollary}  (i) Let $\mu, \nu \in \M_r(G)$. Then $\| \mu * \nu \| = \| \mu \| \| \nu \| $ if and only if $| \mu * \nu | = | \mu|* | \nu|$. \\
(ii) If $\mu, \nu \in \M_r(G) \cap \M(\om)$ and   $\| \mu * \nu \|_\om = \| \mu \|_\om \| \nu \|_\om$, then  $\| \mu * \nu \| = \| \mu \| \| \nu \|$. 
\ec 

\begin{proof}  The forward implication of statement (i)  is a special case of Proposition \ref{Multiplicative Tot Var Prop} (and is also \cite[Theorem 2.1.2]{Gre}). Since $\M_r(G) \ra \C: \mu \mapsto \int 1 \, d \mu_e$ is a homomorphism and $\| \mu \| = |\mu_e|(G) = \int 1 \, d|\mu_e|$, we have the converse of statement (i). Statement (ii) follows immediately from Proposition \ref{Multiplicative Tot Var Prop} and part (i). 
\end{proof}  

\bc \label{Norm-one idempotents in M(w) Corollary} Let $\nu \in \M(\om)$ be a norm-one idempotent. Then there is a compact subgroup $K$ of $G$ and $\rho \in \widehat{K}^1$, i.e., a continuous homomorphism $\rho:K \ra \mathbb{T}$,  such that $\nu = \rho m_K$; moreover $\om \equiv 1$ on $K$. 
 \ec
\begin{proof}  Since $\|\nu * \nu\|_\om = \| \nu\|_\om =1 = \| \nu\|_\om\|\nu\|_\om$, $|\nu| * |\nu| = |\nu*\nu| = |\nu|$ by Proposition \ref{Multiplicative Tot Var Prop}. By Theorem  \ref{Positive Idempotent Thm}, $K = \s(|\nu|) = \s(\nu)$ is a compact subgroup of $G$ and $\om \equiv 1$ on $K$. Hence, $\nu \in \M_{cr}(G) \subseteq \M_r(G)$ and $\|\nu\| =  \| \om \nu\| = \|\nu\|_\om = 1$. Thus, $\nu$ is a norm-one idempotent in $\M_r(G)$, so by \cite[Theorem 2.1.4]{Gre}, $\nu_e = \rho m_K$ --- and therefore $\nu = \rho m_K$ --- for some $\rho \in \widehat{K}^1$. 
\end{proof}

\section{Norm-multiplicative subgroups} 

 To avoid trivial cases, we assume all (convolution) subgroups $\Gamma$ of $\M(\om)$  are not $\{0\}$ (and therefore do not contain 0). Unless stated  otherwise, $\Gamma$ is always endowed with its relative $\sigma^\om$-topology inherited from $\M(\om)$ and   $\iota_\Gamma$ denotes the identity of $\Gamma$.  It is convenient to introduce some terminology.  

\bd \rm We call a subgroup $\Gamma$ of $\M(\om)$  \it norm-multiplicative, \rm or say that $\Gamma$ is a $\NMO$-\it subgroup \rm of $\M(\om)$,  if $$\| \mu * \nu \|_\om = \| \mu \|_\om \| \nu \|_\om \quad  \text{ for } \quad \mu, \nu \in \Gamma,$$i.e., if $\| \cdot \|_\om : \Gamma \ra (0, \infty)$ is a homomorphism;  $\NMO$-subgroups taken with respect to the trivial weight $\om \equiv 1$ will be called $\NM$-subgroups of $\M_r(G)$.  
\ed 

The following observation is an immediate consequence of Corollary \ref{Multiplicative Norm Corollary}.

\bp \label{NM-subgroups Prop} Let  $\Gamma$ be a subgroup of $\M_r(G)$. Then $\Gamma$  is a 
 a $\NM$-subgroup if and only if 
 $| \mu * \nu | = | \mu | * |\nu|$ for every $\mu, \nu \in \Gamma$.
\ep

 In parts (iii) and (iv) of the next proposition, we are not assuming that $\Gamma$ is contained in $\M_r(G)$; in each case, this is a consequence of our hypothesis. 

\bp \label{NMw-subgroup Prop} Let $\Gamma$ be a subgroup of $\M(\om)$. 
\bi \item[(i)]  If $\Gamma$ is a contractive subgroup of $\M(\om)$, i.e., if $\Gamma \subseteq \M(\om)_{\| \cdot\|_\om \leq 1}$, then $\Gamma \subseteq \M(\om)_{\| \cdot\|_\om = 1}$, and $\Gamma$  is therefore a $\NMO$-subgroup of $\M(\om)$.

\item[(ii)]  $\Gamma$ is a $\NMO$-subgroup of $\M(\om)$ if and only if $\displaystyle \Gamma_1=  \left\{ \mu / \|\mu\|_\om: \mu \in \Gamma\right\}$ is a contractive subgroup of $M(\om)$.  

\item[(iii)] If $\Gamma$ is a $\NMO$-subgroup of $\M(\om)$, then $\|\ig\|_\om =1$,  $\Gamma \subseteq \M_{cr}(G)$ and $\Gamma$ is a $\NM$-subgroup of $\M_r(G)$. 

\item[(iv)] If $\Gamma \subseteq \M(\om)^+$ and $\ig \in \M_r(G)$ --- e.g., if $\|\ig\|_\om =1$ --- then $\Gamma$ is a $\NM$-subgroup of $\M_r(G) $ and $\iota_\Gamma = m_K$ for some compact subgroup $K$ of $G$. 
\ei
\ep 

\begin{proof} (i) As a nonzero idempotent in $\M(\om)$, $\|\ig\|_\om \geq 1$, so $\|\ig\|_\om =1$. For $\mu \in \Gamma$, it follows that $1 = \|\ig\|_\om = \|\mu * \mu^{-1}\|_\om  \leq \|\mu \|_\om \| \mu^{-1}\|_\om \leq 1$; hence $\| \mu \|_\om =1$. 

\smallskip 

\noindent (ii) Since $\Gamma_1$ is contained in the unit sphere  of $M(\om)$,  this is easy to verify.  

\smallskip 

 \noindent (iii) Clearly $\| \ig\|_\om =1$, so by Corollary \ref{Norm-one idempotents in M(w) Corollary}, $\s(\ig) =K$ is a compact subgroup of $G$. Let $\mu \in \Gamma$ and take $y \in \s(\mu^{-1})$. By Proposition \ref{Multiplicative Tot Var Prop}, $K = \s(\ig) = \ov{\s(\mu) \s(\mu^{-1})}$, so $xy \in K$  for $x \in \s(\mu)$ and therefore $x \in Ky^{-1}$. Hence, $\mu \in \M_{cr}(G) \subseteq \M_r(G)$. Statement (i) now follows from Corollary \ref{Multiplicative Norm Corollary}. 

\smallskip

\noindent (iv) Since $\ig \in \M_r(G)^+ \cong M(G)^+$, $\| \ig \| = \l 1_G, \ig \r = \l 1_G, \ig * \ig \r = \l 1_G, \ig \r^2  = \|\ig\|$. Hence, $\ig$ is a positive norm-one idempotent in $\M_r(G)$, so $K = \s(\ig)$ is a compact subgroup of $G$ and $\ig = m_K$ by the unweighted case of  Theorem \ref{Positive Idempotent Thm} (or by \cite[Theorem 4.1]{Pym}). By \cite[Proposition III.10.5]{Fel-Dor}, $\s(\mu * \nu) = \ov{\s(\mu)\s(\nu)}$ for $\mu, \nu \in \Gamma$ so the argument given in (iii) above shows that $\Gamma \subseteq \M_{cr}(G)$. By Proposition \ref{NM-subgroups Prop}, $\Gamma$ is a $\NM$-subgroup of $\M_r(G)$.
\end{proof} 

\br \rm Observe that  $\Gamma = \{ \delta_x: x \in G\}$ is a positive subgroup of $\M(\om)$ with $\| \ig\|_\om =1$, and is therefore a $\NM$-subgroup of $\M_r(G)$. However, it is easy to see that $\Gamma$ is a $\NMO$-subgroup of $M(\om)$ if and only if $\om$ is multiplicative on $G$; more generally, see Theorem \ref{NMO-subgroup Description Theorem} below.  We do not know if a positive idempotent in $\M(\om)$  must always lie in $\M_r(G)$, though when $\om$ is bounded away from zero, $\M(\om)$ is  automatically contained in $ \M_r(G)$. 
\er

We now observe that Proposition \ref{NMw-subgroup Prop} and the description of the contractive subgroups of $M(G) \cong \M_r(G)$ from \cite{Gre} and \cite{Sto1} can be used to describe the $\NM$-subgroups of $\M_r(G)$. Our notation follows \cite[Section 4]{Sto1}, which is slightly different from Greenleaf's. 

Let $\Gamma$ be a $\NM$-subgroup of $\M_r(G)$. Since $\| \iota_\Gamma\| = 1$, there is a compact subgroup $K$ of $G$ and $\rho \in \widehat{K}^1$ such that $\iota_\Gamma = \rho m_K$  and $\rho m_K$ is self-adjoint \cite[Lemma 4.1]{Sto1}. By Proposition \ref{NMw-subgroup Prop},
$$\Gamma_1=  \left\{ \mu_1:= \mu / \|\mu\|: \mu \in \Gamma\right\}$$
is a contractive subgroup of $\M_r(G)$.  From Greenleaf's description of such groups \cite[Section 3]{Gre}: 
\bi \item $ \ds H_0 = \bigcup_{\mu \in \Gamma} \s(\mu) = \bigcup_{\mu \in \Gamma} \s(\mu_1)$ is a subgroup of $G$ and, letting $H = \s(\Gamma):=\ov{H_0}$, $K$ and $\ker \rho$ are normal subgroups of $H$ such that $K/\ker \rho$ is contained in $Z(H/\ker \rho)$, the centre of $H/\ker \rho$. As noted in \cite[Lemma 5.2]{Sto1}, $\iota_\Gamma \in Z(\M_r(H))$, the centre of $\M_r(H)$ (and $\M_r(H)$ is a closed subalgebra of  $ \M_r(G)$, via Remark \ref{Remark mu_e}.2). 

\item Letting $$\Omega_{\Gamma_1} = \{(\alpha, t) \in \T \times H: \alpha \delta_t * \rho m_K \in \Gamma_1\} \quad \text{and} \quad \Omega_\rho= \{(\rho(k),k): k \in K \}, $$
$\Omega_{\Gamma_1}$ is a subgroup of $\T \times H$,  and $\Gamma_{\T \times H}= \{ \alpha \delta_t * \rho m_K: (\alpha, t) \in \T \times H\}$ is a contractive subgroup of $\M_r(H) \subseteq \M_r(G)$ containing $\Gamma_1$. 
\item The map $ \phi_0^1: (\alpha, t) \mapsto \alpha \delta_t * \rho m_K$ 
defines a continuous homomorphism of $\Omega_{\Gamma_1}$ onto $\Gamma_1$, and of $\T \times H$ onto  $\Gamma_{\T \times H} $,  with $\ker\phi^1_0 = \Omega_\rho$; 
$$\phi^1: (\alpha, t) \Omega_\rho \mapsto \alpha \delta_t * \rho m_K$$ is thus a continuous group isomorphism of $\Omega_{\Gamma_1}/\Omega_\rho$ onto $\Gamma_1$, and of $\T \times H /\Omega_\rho$ onto $\Gamma_{\T \times H}$. Moreover, $\phi^1$ is a topological group isomorphism \cite[Theorem 4.2]{Sto1}, \cite[Theorem 3.1]{Spr}.
\ei 

\br \rm \label{Properties of alpha deltas * rho mK Remark}  1. (a)  For $s, t \in H$ and $\alpha, \beta \in \C$, 
$(\alpha \delta_s * \rho m_K) * (\beta \delta_t * \rho m_K) = \alpha\beta \,  \delta_{st} * \rho m_K$
because $\iota_\Gamma = \rho m_K \in Z(\M_r(H))$.  \\
(b) For $\mu = \alpha \delta_s * \rho m_K$ with $\alpha \in \C$ and $s \in G$, 
$\| \mu \| = | \alpha| \| \delta_s * \rho m_K \| = | \alpha|$; for any $\mu = \delta_{s^{-1}} * \delta_s * \mu \in \M_r(G)$, $\| \delta_s * \mu \| = \| \mu\| $ by submultiplicativity of $\| \cdot \|$.  \\
(c) From (b), if $ \alpha \delta_s * \rho m_K = \rho m_K$, then $\alpha \in \T$ and therefore $(\alpha, s) \in \ker \phi^1_0 = \Omega_\rho$.  

\smallskip 

\noindent 2. Let $\mu = \alpha \delta_s * \rho m_K$  for $\alpha \in \C$ and $s \in G$, where we further assume that $\om \equiv 1$ on $K$.  Then $\mu \in \M_{cr}(G) \subseteq \M(\om)$ and for $\phi \in C_0(G)$,
\beqs \l \om \mu, \phi \r  = \int_K \alpha \om(sk) \phi(sk) \rho(k) \, dm_K(k) =  \om(s) \int_K \alpha \phi(sk) \rho(k) \, dm_K(k) = \l \om(s) \mu, \phi \r, \eeqs 
using Remark \ref{Constant on Cosets Remark Etc}. Hence, in this case we have
\beq \label{omega mu eqn} \om \mu = \om(s) \mu \quad \text{and} \quad \|\mu \|_\om = \| \om \mu \| = \om(s) \| \mu \| = \om(s)|\alpha|. \eeq
\er 

Let  $\Cx$ denote the group $\C \bs \{0\}$ under multiplication. 

\blem \label{Main NM-subgroup Description Lemma} Let $\Gamma$ be a $\NM$-subgroup of $\M_r(G)$. Then: 
\bi \item[(i)] there is a compact subgroup $K$ of $G$ and $\rho \in \widehat{K}^1$ such that $\iota_\Gamma = \rho m_K$, $K$ and $\ker \rho$ are normal subgroups of $H$ such that $K / \ker \rho \leq Z(H/\ker \rho)$ and $\rho m_K$ is a self-adjoint central idempotent in $\M_r(H) \subseteq \M_r(G)$; 

\item[(ii)] $\Omega_\Gamma = \{(\alpha, t) \in \Cx \times H: \alpha \delta_t * \rho m_K \in \Gamma\}$    is a subgroup of   $\Cx \times H$ and 
\beq \label{phi_0 Def Eqn}  \phi_0: (\alpha, t) \mapsto \alpha \delta_t * \rho m_K\eeq
defines a continuous homomorphism of $\Omega_\Gamma$ onto $\Gamma$ with $\ker\phi_0 = \Omega_\rho$ and $\mu^{-1} =\mu^* / \|\mu\|^2$ for $\mu \in \Gamma$;

\item[(iii)]  $\Gamma_{\Cx \times H}= \{ \alpha \delta_t * \rho m_K: (\alpha, t) \in \Cx \times H\}$ is a $\NM$-subgroup of $\M_r(G)$ and   (\ref{phi_0 Def Eqn}) also defines a continuous homomorphism of $\Cx \times H$ onto $\Gamma_{\Cx \times H}$ with $\ker \phi_0 = \Omega_\rho$.

\ei 
 \elem

\begin{proof}  Statement (i) was already observed,  and the algebraic parts of statement (iii) follow from Remark \ref{Properties of alpha deltas * rho mK Remark}.1 (a), (b) and (c). Since $s \mapsto \delta_s$ is a topological group isomorphism of $H$ onto $(\Delta_H, \sigma)$ and multiplication in $\M_r(H)$ is separately $\sigma$-continuous,  $\phi_0$ is continuous on $\Cx\times H$. That $\Omega_\Gamma$ is a subgroup of $\Cx \times H$ and $\phi_0$ is a continuous homomorphism of $\Omega_\Gamma$ into $\Gamma$, with kernel $\Omega_\rho$, now follows from Remark \ref{Properties of alpha deltas * rho mK Remark}.1 (a), (c) and statement (iii). If $\mu \in \Gamma$, $\mu_1 = \mu / \|\mu\| \in \Gamma_1$, so $\mu_1 = \alpha \delta_s * \rho m_K$ for some $(\alpha, s ) \in \Omega_{\Gamma_1}$; therefore, $(\|\mu\| \alpha, s) \in \Omega_\Gamma$ and $\phi_0(\|\mu\| \alpha, s) = \mu$. 
\end{proof}

\bt \label{NMO-subgroup Description Theorem}   Let $\Gamma$ be a subset of $\M_r(G)$. The following three statements are equivalent: 
\bi \item[(i)]  $\Gamma$  is a $\NMO$-subgroup of $\M(\om)$.
\item[(ii)] $H = \s(\Gamma)$  is a closed subgroup of $G$ on which $\om: H \ra (0,\infty)$ is a continuous homomorphism; there is a compact normal subgroup $K$ of $H$ and $\rho \in \widehat{K}^1$ such that $\ker \rho \lhd H$ and $K/\ker \rho \leq Z(K /\ker \rho)$; and $\Omega_\Gamma = \{(\alpha, s) \in \Cx \times H: \alpha \delta_s * \rho m_K \in \Gamma\}$ is a subgroup of $\Cx \times H$ such that $\Gamma = \{ \alpha \delta_s * \rho m_K: (\alpha, s) \in \Omega_\Gamma\}$.
\item[(iii)] $\Gamma $ is a $NM$-subgroup of $\M_r(G)$ and $\om$ is a continuous homomorphism on $H = \s(\Gamma)$. 
\ei
Hence,  $\Gamma_{\Cx \times H}= \{ \alpha \delta_s * \rho m_K: (\alpha, s) \in \Cx \times H\}$ is a $\NMO$-subgroup of $\M(\omega)$ whenever $\Gamma$  is a $\NMO$-subgroup of $\M(\om)$ with $H = \s(\Gamma)$.
\et 

\begin{proof}  (i) $\Rightarrow$ (ii): By Proposition \ref{NMw-subgroup Prop} (iii), $\Gamma$ is a $\NM$-subgroup of $\M_r(G)$ and $\|\iota_\Gamma\|_\om =1$, so by Lemma \ref{Main NM-subgroup Description Lemma} we only need to show that $\om$ is multiplicative on the dense subgroup $H_0 = \cup_{\mu \in \Gamma} \, \s(\mu)$ of $H$. Let  $s, t \in H_0$ and, again using Lemma \ref{Main NM-subgroup Description Lemma},  observe that there must exist some $\alpha, \beta \in \Cx$ such that $\mu = \alpha \delta_s * \rho m_K, \nu = \beta \delta_t * \rho m_K \in \Gamma$. By Corollary \ref{Norm-one idempotents in M(w) Corollary}, $\om \equiv 1$ on $K$, so (\ref{omega mu eqn}) yields 
\beqs \om(st) \| \mu \| \| \nu \| &=&  \om(st) \| \mu * \nu \| = \om(st) \| \alpha \beta \delta_{st} * \rho m_K\| = \| \om (\mu * \nu) \| \\ & = &  \|\mu * \nu \|_\om  =  \| \mu \|_\om \| \nu \|_\om = \om(s) \| \mu \| \om(t) \|\nu\|.\eeqs
We conclude that  $\om(st) = \om(s) \om(t)$. 

\smallskip 

\noindent (ii) $\Rightarrow$ (iii): Observe that $\Gamma \subseteq \M_{cr}(H)$ and $\rho m_K$ is a central idempotent in $\M_r(H)$ by \cite[Lemma 5.2]{Sto1}. Hence, Remark \ref{Properties of alpha deltas * rho mK Remark}.1 (a) holds for all $(\alpha, s), (\beta, t) \in \Omega$, and it follows that $\Gamma$ is a subgroup of $\M_r(G)$ with identity $\iota_\Gamma = \rho m_K$. Moreover, by Remark \ref{Properties of alpha deltas * rho mK Remark}.1 (b), $\Gamma$ is a $\NM$-subgroup of $\M_r(G)$. 

\smallskip 

\noindent (iii) $\Rightarrow$ (i): Let $\mu, \nu \in \Gamma$, say $\mu = \alpha \delta_s * \rho m_K$ and $\nu = \beta \delta_t * \rho m_K$, where we are using Lemma \ref{Main NM-subgroup Description Lemma}. Since $\om$ is a continuous homomorphism on $H$, $\om \equiv 1$ on the compact subgroup $K$, and (\ref{omega mu eqn})  yields 
\beqs \|\mu * \nu \|_\om & =&  \| \om( \alpha \beta \delta_{st}* \rho m_K)\|  = \om(st) \| \alpha \beta \delta_{st} * \rho m_K \| \\ & =&  \om(st) \| \mu * \nu \| = \om(s) \om(t) \| \mu \| \| \nu\| = \| \mu\|_\om \| \nu \|_\om. \qedhere \eeqs
\end{proof} 

As noted in Section 1, the $so_l^\om$-,  $so_r^\om$-,  $so^\om$- and weak$^*$-topologies agree on $\Delta_G \subseteq \M(\om)$. Statement (v) of Proposition \ref{NM-subgroup Isomorphism Prop}, show that this   is a property of $\NM$-subgroups in general.

\bp \label{NM-subgroup Isomorphism Prop}  Let $\Gamma$  be a $\NM$-subgroup of $\M_r(G)$, with $$H = \s(\Gamma),  \quad \iota_\Gamma = \rho m_K \quad \text{and} \quad  \Omega_\Gamma = \{ (\alpha, s) \in \Cx \times H : \alpha \delta_s * \rho m_K\in \Gamma\}$$  as described in Lemma \ref{Main NM-subgroup Description Lemma}.  Let $\om$ be any weight on $G$.  Then: 

\bi \item[(i)] $\Gamma, \ \Gamma_{\Cx \times H} \subseteq \M_{cr}(G) \subseteq \M(\om) \cap \M_r(G)$  are subgroups of $\M(\om)$ and $\phi_0: (\alpha, s) \mapsto \alpha \delta_s * \rho m_K$ is a continuous homomorphism mapping $\Omega_\Gamma$ onto $(\Gamma, \sigma^\om)$ and $\Cx \times H$ onto $(\Gamma_{\Cx \times H}, \sigma^\om)$, with kernel $\Omega_\rho = \{ (\rho(k),k): k \in K \}$. 

\item[(ii)] $\phi: (\alpha, s)\Omega_\rho \mapsto \alpha \delta_s * \rho m_K$ defines a topological group isomorphism of $\Omega_\Gamma/ \Omega_\rho$ onto $(\Gamma, \sigma^\om)$ and $\Cx \times H/\Omega_\rho$ onto $(\Gamma_{\Cx \times H}, \sigma^\om)$. 

\item[(iii)] The  relative $\sigma$ and $\sigma^\om$-topologies all agree on $\Gamma$ and $\Gamma_{\Cx \times H}$; so, we can unambiguously refer to the weak$^*$-topology on a $\NM$-subgroup of $\M_r(G)$, (or $\NMO$-subgroup of $\M(\om)$).

\item[(iv)]  $\Gamma$ is a topological group and $\Gamma_{\Cx \times H} $ is a locally compact group.

\item[(v)] The $so_l^\om$-,  $so_r^\om$-,  $so^\om$- and weak$^*$-topologies all agree on $\Gamma$. 

\ei 

\ep

\begin{proof}  (i) This mostly follows from Lemma \ref{Main NM-subgroup Description Lemma}, with $\sigma^\om$-continuity of $\phi_0$ following from the $\sigma^\om$-continuity of $s \mapsto \delta_s$ and separate $\sigma^\om$-continuity of multiplication in the dual Banach algebra $\M(\om)$. 

\smallskip 

\noindent (ii) From (i), $\phi$ is a continuous group isomorphism. Moreover, the proof of \cite[Theorem 4.2]{Sto1}, verbatim, shows that $\phi^{-1}: (\Gamma, \sigma^\om) \ra \Omega_\Gamma/\Omega_\rho$ is continuous,  because the ``$v$" from that proof belongs to $C_{00}(H)$ (and this $v$ can be extended to $v_0\in C_{00}(G) \subseteq C_0(\om^{-1})$ if one prefers to work on $(G, \om)$ rather than $(H, \om{\vert}_H)$). 

\smallskip 

\noindent (iii) and (iv) are an immediate consequence of (ii). 

\smallskip

\noindent (v) We will show that the $so_l^\om$- and $\sigma^\om$-topologies agree on $\Gamma$. By \cite[Proposition 2.9]{Sto2}, the identity $\M(\om) \ra \M(\om)$ is $so_l^\om$-$\sigma^\om$ continuous --- take ${\cal S}(\om^{-1}) = C_0(\om^{-1})$ there --- so it suffices to establish continuity of $\text{id}: (\Gamma, \sigma^\om) \ra (\Gamma, so_l^\om)$. Since $(\Gamma, \sigma^\om)$ is a topological group and multiplication in $\M(\om)$ is separately $so_l^\om$-continuous, it suffices to establish continuity at $\iota_\Gamma = \rho m_K$. To this end, let $(\mu_i)$ be a net in $\Gamma$ such that $\mu_i \ra \rho m_K$ $\sigma^\om$, say $\mu_i = \alpha_i \delta_{s_i} * \rho m_K$ with $(\alpha_i, s_i) \in \Omega_\Gamma \subseteq \Cx \times H$. 
Then $\mu_i \ra \rho m_K$ in $\Gamma_{\Cx \times H}$ as well, so $(\alpha_i, s_i) \Omega_\rho = \phi^{-1}(\mu_i) \ra  \phi^{-1}(\rho m_K) = \Omega_\rho$ in $\Cx \times H/\Omega_\rho$. Take a relatively compact neighbourhood $W$ of $\Omega_\rho$ in $\Cx \times H$ such that $(\alpha_i, s_i)$ is eventually in $W$;  e.g., letting $W_0$ be any relatively compact neighbourhood of $\Omega_\rho$ in $\Cx \times H$, take $W = W_0 \Omega_\rho \subseteq \Cx \times H$. Let $((\alpha_{i_j}, s_{i_j}))_j$ be a subnet and $(\alpha, s) \in \overline{W}$ such that $(\alpha_{i_j}, s_{i_j}) \ra (\alpha, s)$ in $\Cx \times H$. Then $(\alpha_{i_j}, s_{i_j})\Omega_\rho \ra (\alpha, s) \Omega_\rho$ and  $(\alpha_{i_j}, s_{i_j})\Omega_\rho \ra \Omega_\rho,$
so $(\alpha, s) = (\rho(k_0), k_0)$ for some $k_0 \in K$. As noted in the introduction, this means that $\delta_{s_{i_j}} \ra \delta_{k_0}$ $so_l^\om$ and therefore, again using separate $so_l^\om$-continuity of  multiplication in $\M(\om)$, $\mu_{i_j} = \alpha_{i_j} \delta_{s_{i_j}}* \rho m_K  \ra \rho(k_0) \delta_{k_0}* \rho m_K = \rho m_K \ so_l^\om$ in $\M(\om)$. \end{proof}

\br \rm 1. Spronk \cite[Theorem 3.1(ii)]{Spr} establishes continuity of $\phi^{-1} : \Gamma_{\T \times H} \ra \T \times H/\Omega_\rho$ using a different argument from the one in \cite{Sto1}. This argument seems to require prior knowledge that the scalars ``$z_i$" from that proof are bounded and  does not seem to adapt to our situation.  

\noindent 2.  If we had assumed that $\| \cdot \|: \Gamma \ra (0, \infty)$ is continuous in our definition of $\NM$-groups, we could establish that $\phi^{-1}$ in Proposition \ref{NM-subgroup Isomorphism Prop}(ii)  is continuous by borrowing directly from   the statement of \cite[Theorem 4.2]{Sto1} alone.   In fact,   continuity of $\| \cdot \|_\om$ is automatic on any $\NMO$-subgroup of $\M(\om)$.    
\er 

\bc Let $\Gamma$ be a $\NMO$-subgroup of $\M(\om)$. Then $\| \cdot \|_\om: \Gamma \ra (0,\infty)$ is a  continuous homomorphism. \ec 

\begin{proof}  First,  we \it claim \rm that $\| \cdot \| : \Gamma \ra (0, \infty)$ is continuous. To this end, suppose that $\mu_i \ra \mu$ in $\Gamma$ and therefore  $\mu_i \ra \mu$ in $\Gamma_{\Cx \times H}$, a $\NM$-subgroup of $\M_r(G)$. Observe that for $\nu = \alpha \delta_s * \rho m_K$, $\nu * \nu^* = |\alpha|^2 \rho m_K = \| \nu\|^2 \rho m_K$. Since $\mu_i^*, \mu^* $ also belong to $\Gamma_{\Cx \times H}$, a topological group with continuous involution, 
$$\| \mu_i\|^2 \rho m_K = \mu_i * \mu_i^* \ra \mu * \mu^* = \| \mu \|^2 \rho m_K$$
in $\Gamma_{\Cx \times H}$, and therefore weak$^*$ in $\M_r(G)$. Taking $v \in C_{00}(G) $ such that $v{\vert}_K = \overline{\rho}$, we obtain  
$$\| \mu_i \|^2 = \l \| \mu_i \|^2  \rho m_K, v\r \ra \l \| \mu \|^2 \rho m_K, v \r = \| \mu\|^2,$$
establishing the claim.  Now observe that $\om_{H/K} : H/K \ra (0,\infty)$ is a well-defined continuous homomorphism, as is $p: \Cx \times H / \Omega_\rho \ra H/K: (\alpha, s) \Omega_\rho \mapsto sK$, and $\| \mu\|_\om = \om_{H/K}(p(\vp^{-1}(\mu)))\| \mu \|$ holds for $\mu \in \Gamma$ by Remark \ref{Properties of alpha deltas * rho mK Remark}.2.
\end{proof} 

Let $\R^+ = (0, \infty)$,  the multiplicative group of positive real numbers.  The next corollary contains \cite[Proposition 1.3(i)]{Gha-Zad}. (We remark that   \cite[Proposition 1.3(i)]{Gha-Zad}  also quickly follows   from the version of \cite[Proposition III.10.5]{Fel-Dor} for positive measures.)    

\bc (i) Let $\Gamma$ be a subgroup of $\M(\om)^+$, the set of positive measures in $\M(\om)$, with $\iota_\Gamma \in \M_r(G)$.  Then there is a compact normal subgroup $K$ of $H = \s(\Gamma)$ such that  $\Omega_\Gamma = \{ (\alpha, s)\in \R^+ \times H : \alpha \delta_s * m_K\in \Gamma \} $ is a subgroup of $\R^+ \times H$ such that $\phi: (\alpha, s) (1 \times K) \mapsto \alpha \delta_s * m_K$ is topological group isomorphism of $\Omega_\Gamma /(1 \times K)$ onto $\Gamma$. Moreover, $\Gamma_{\R^+ \times H} = \{ \alpha \delta_s * m_K: (\alpha , s ) \in \R^+ \times H\}$ is a subgroup of $\M(\om)^+$ that is isomorphic as a topological group  with $\R^+ \times (H /K)$. 

\smallskip 

\noindent (ii) If $\Gamma $ is any subgroup of invertible elements in $\M(\om)^+$, then $\Omega_\Gamma = \{ (\alpha, s)\in \R^+ \times H : \alpha \delta_s \in \Gamma \} $ is a subgroup of $\R^+ \times G$ such that $\phi: (\alpha, s)  \mapsto \alpha \delta_s $ is topological group isomorphism of $\Omega_\Gamma $ onto $\Gamma$. In particular, if ${\rm PInv}(G,\om)$ is the group of all positive invertible elements in $\M(\om)$ with positive inverses, then ${\rm PInv}(G,\om) = \{ \alpha \delta_s : (\alpha, s) \in \R^+ \times G\} \cong \R^+ \times G$.
\ec

\begin{proof}  Let $\Gamma$ be a subgroup of $\M(\om)^+$ with $\iota_\Gamma \in \M_r(G)$. By Proposition \ref{NMw-subgroup Prop} (iv), $\Gamma$ is a $\NM$-subgroup  of $\M_r(G)$ with $\iota_\Gamma = m_K$, so, $\Omega_\rho = 1\times K$ and $\Omega_\Gamma \leq \R^+ \times H$, so the first part of statement (i) follows from  Proposition \ref{NM-subgroup Isomorphism Prop}. Since $\Gamma_{\R^+ \times H} \leq \Gamma_{\Cx \times H}  \cong (\Cx \times H)/(1 \times K) \cong \Cx \times H/K$, we obtain the second part of statement (i). Part (ii) is an immediate consequence of part (i) because $K = \{e_G\}$ in this case, and $\{ \alpha \delta_s : (\alpha, s) \in \R^+ \times G\}$ is contained in  ${\rm PInv}(G,\om)$. 
\end{proof} 


\section{Norm-multiplicative homomorphisms} 

Throughout this section, $(F, \omf)$ and $(G, \omg)$ are weighted locally compact groups. 

\br \label{Canonical extension to M(omf) Remark}  \rm By \cite[Corollary 2.12]{Sto2}, every bounded homomorphism $\vp: \L^1(\omf) \ra \M(\omg)$ has a unique $so_l^\omf-\sigma^\omg$  continuous --- equivalently $so_r^\omf-\sigma^\omg$ continuous \cite[Corollary 2.13]{Sto2}  --- extension to  $\M(\omf)$, also denoted herein by $\vp$, with the same norm. Thus, the bounded homomorphisms $\vp: \L^1(\omf) \ra \M(\omg)$ and $so_l^\omf-\sigma^\omg$ continuous homomorphisms  $\vp: \M(\omf) \ra \M(\omg)$ are in one-to-one correspondence, and our theorems describing certain homomorphisms $\vp:\L^1(\omf) \ra \M(\omg)$ equivalently describe $so_l^\omf-\sigma^\omg$ continuous homomorphisms $\vp: \M(\omf) \ra \M(\omg)$. If $\mu \in \M(\omf)$ and  $(f_i)$ is a bounded approximate identity for $\L^1(\omf)$, $(f_i * \mu)$ is a net in $\L^1(\omf)$ such that 
\beq \label{so_l-wk* conts extension identity} so_l^\omf - \lim f_i * \mu = \mu, \qquad \text{so} \qquad \vp (\mu) = \sigma^\omg-\lim \vp(f_i*\mu). \eeq
Taking $(f_i)$ to be a positive contractive approximate identity for $\L^1(\omf)$ \cite[Proposition 1.6]{Sto2}, it follows that $\vp$ is positive or contractive on $\M(\omf)$ when the same is true for $\vp$ on $\L^1(\omf)$. 
\er 

\bd \label{NM-Homoms Defn}  \rm A  bounded $so_l^\omf-\sigma^\omg$ continuous homomorphism $\vp:\M(\omf) \ra \M(\omg)$, or a  bounded  homomorphism $\vp: \L^1(\omf) \ra \M(\omg)$  with unique $so_l^\omf-\sigma^\omg$ continuous homomorphic extension $\vp: \M(\omf) \ra \M(\omg)$, will be called a \it  $\NM_\omg$-homomorphism \rm if $\Gamma_\vp:= \vp(\Delta_F)$ is a $\NM_\omg$-subgroup of  $\M(\omg)$;  $\vp$ is a \it $\NM$-homomorphism \rm if $\Gamma_\vp$ is a $\NM$-subgroup of $\M_r(G)$. 
\ed 

\br \label{NM vs NM_omega Homoms Remark} \rm 1.  We will mostly focus on $\NM$-homomorphisms and \it stress \rm that our definition of a $\NM$-homomorphism does not assume $\vp$ maps  $\L^1(\omf)$ into  $\M_r(G)$; we assume that $\vp: \L^1(\omf) \ra \M(\omg)$, but that $\Gamma_\vp \subseteq \M_r(G)$ (trivially satisfied when $\omg$ is bounded away from zero). 

\smallskip 

\noindent 2.  By   Proposition \ref{NMw-subgroup Prop}, a homomorphism $\vp: \L^1(\omf) \ra \M(\omg)$ is a $\NM$-homomorphism in any of the following cases: \bi 
\item[(i)] $\vp$ is a $\NM_\omg$-homomorphism (in which case $\| \vp (\delta_{e_F})\|_\omg=1$);
\item[(ii)] $\omf\equiv 1$ and $\vp$ is a contraction (a special case of (i));
\item[(iii)] $\Gamma_\vp \subseteq \M(\omg)^+$ and $\vp(\delta_{e_F}) \in \M_r(G)$ (e.g, if $\| \vp(\delta_{e_F}) \|_\omf =1$). 
\ei 
Thus, our characterization below of $\NM$-homomorphisms includes all contractive homomorphisms $\vp: \L^1(F) \ra \M(\omg)$ --- described in the $\omg\equiv 1$ case in \cite{Gre} and \cite{Sto1} (papers on which the present work depends) --- and \it all \rm bounded positive homomorphisms $\vp: \L^1(F) \ra \M_r(G)$, a new result. 
\er  
We do not know if compositions of $\NM$-homomorphisms are always  $\NM$-homomorphisms, but do have the following proposition (where $(Q, \om_Q)$ is another weighted locally compact group). Note that a $\NM$-subgroup of $\M_r(F)$ is automatically contained in $\M_{cr}(F) \subseteq \M(\omf)$ by Lemma \ref{Main NM-subgroup Description Lemma}. 

\bp \label{Composition of NM-homs Prop} Let $\vp: \M(\omf) \ra \M(\omg)$ and $\kappa: \M(\omg) \ra \M(\om_Q)$ be $\NM$-homomorphisms. 
\bi \item[(i)]  If  $\Gamma$ is a $\NM$-subgroup of $\M_r(F)$ and $\vp(\iota_\Gamma) \in \M_r(G)_{\| \cdot\|\leq1}$  --- e.g., if  $\|\vp(\iota_\Gamma)\|_\omg \leq 1$ --- then $\vp(\Gamma)$ is a $\NM$-subgroup of $\M_r(G)$. 

\item[(ii)] If  $\kappa$ is weak$^*$-continuous and  $\kappa \circ \vp (\delta_{e_F}) \in \M_r(Q)_{\| \cdot\|\leq1} $, then $\kappa \circ \vp$ is a $\NM$-homomorphism. 
\ei 
\ep  

\begin{proof}  Let $\Lambda = \vp (\Gamma)$, a subgroup of $\M(\omg)$, with identity $\iota_\Lambda = \vp(\iota_\Gamma)$. Since $\iota_\Lambda$ is an idempotent and $\| \iota_\Lambda \| =1$ by assumption,  $\iota_\Lambda \in \M_{cr}(G)$. Let $\nu_1, \nu_2$ belong to the $\NM$-subgroup $\Gamma$. By Proposition \ref{NM-subgroup Isomorphism Prop}, for $i =1,2$ we have $\nu_i  = \alpha_i \delta_{s_i} * \iota_\Gamma$ for some $(\alpha_i, s_i) \in \Cx \times \s(\Gamma)$; also $\vp(\delta_{s_i}) = \beta_i \delta_{t_i}* \iota_\Phi$ for some $(\beta_i, t_i) \in \Cx \times \s(\Phi)$, since $\Phi = \vp(\Delta_F)$ is a $\NM$-subgroup of $\M_r(G)$. Observe that $\iota_\Phi * \iota_\Lambda = \vp(\delta_{e_F}) * \vp(\iota_\Gamma) = \vp(\delta_{e_F} * \iota_\Gamma) = \iota_\Lambda$, so $\vp(\nu_i) = \alpha_i \vp(\delta_{s_i}) * \vp(\iota_\Gamma) = \alpha_i \beta_i  \delta_{t_i} * \iota_\Lambda$. Hence, $\Lambda$ is contained in $\M_{cr}(G) \subseteq \M_r(G)$ and, using properties of multiplication in the $\NM$-subgroups $\Gamma$ and $\Phi$, 
\beqs  \vp(\nu_1) * \vp(\nu_2) & = &  \vp (\nu_1*\nu_2) = \vp( \alpha_1 \alpha_2 \,  \delta_{s_1} * \delta_{s_2} * \iota_\Gamma) = \alpha_1\alpha_2\, \vp(\delta_{s_1}) * \vp(\delta_{s_2}) * \vp(\iota_\Gamma)   \\  &=&  \alpha_1\alpha_2 (\beta_1 \delta_{t_1}* \iota_\Phi * \beta_2 \delta_{t_2}* \iota_\Phi)*\iota_\Lambda = \alpha_1\alpha_2 \beta_1 \beta_2\,  \delta_{t_1t_2}* \iota_\Phi*\iota_\Lambda \\ & = & \alpha_1\beta_1\alpha_2  \beta_2 \, \delta_{t_1t_2}* \iota_\Lambda.  \eeqs
By Remark \ref{Properties of alpha deltas * rho mK Remark}.1(b), $\|  \vp(\nu_1) * \vp(\nu_2)\| =\| \vp(\nu_1)\|\| \vp(\nu_2)\| $. This establishes statement (i), from which statement (ii) immediately follows. 
\end{proof}  

\subsection{\rm The basic NM-homomorphisms $j_{\gamma, \theta}^*$, $j_\theta^*$, $M_\gamma^*$   and $S_K^*$} 

Let $\gamma: F \ra \Cx$ and  $\theta: F \ra G$
be continuous homomorphisms and put 
$$\Jgt(\psi) = \gamma(\psi \circ \theta) \qquad \text{for}  \ \psi \in LUC(\omg^{-1}).$$

\bp  The linear map $\Jgt: LUC(\omg^{-1}) \ra LUC(\omf^{-1})$ is well-defined and bounded if and only if $$L_{\gamma, \theta} := \sup_{x \in F} \left|\gamma(x)\right| {\omg(\theta(x)) \over \omf(x)} < \infty.$$ In this case, $\|\Jgt\| = L_{\gamma, \theta} $ and $\Jgt^*: LUC(\omf^{-1})^* \ra LUC(\omg^{-1})^*$ is the unique weak$^*$-continuous homomorphism mapping $\delta_x$ to  $\gamma(x) \delta_{\theta(x)}$ ($x \in F$). 
\ep 

\begin{proof}  Observe that $\omg \in LUC(\omg^{-1})$, and  $\Jgt(\omg) \in \ell^\infty(\omf^{-1})$ if and only if $L_{\gamma, \theta}  <\infty$. Suppose that $L_{\gamma, \theta}  <\infty$. Then 
$$\left|{\Jgt \psi(x) \over \omf(x)} \right| = \left|\gamma(x)\right| {\omg(\theta(x)) \over \omf(x)}    \left| { \psi( \theta(x)) \over \omg(\theta(x))} \right|  \leq L_{\gamma, \theta}  \| \psi\|_{\infty, \omg^{-1}}$$
for $\psi \in LUC(\omg^{-1})$, $x \in F$; hence  $\Jgt: LUC(\omg^{-1}) \ra \ell^\infty (\omf^{-1})$ and  $\|\Jgt \| = L_{\gamma,\theta}$.  For $s \in F$ and $\psi \in LUC(\omg^{-1})$, $(\Jgt \psi) \c s = \gamma(s) \Jgt(\psi \c \theta(s))$ and $\| \psi \c \theta(s)\|_{\infty, \omg^{-1}} \leq \|\psi \|_{\infty, \omg^{-1}} \omg(\theta(s)) $, from which one can check that $\Jgt(\psi) \in LUC(\omf^{-1})$ using \cite[Lemma 2.3(b)]{Sto2}. It is easy to check that $\Jgt^*(\delta_x) = \gamma(x) \delta_{\theta(x)}$, so  $\Jgt^*$ is a homomorphism by Lemma \ref{WK* continuous Homom Lemma}.
\end{proof}

 Let   $\jgt= \Jgt\vert_{C_0(\omg^{-1})} = \Jgt \circ \iota_{C_0}$, where $\iota_{C_0}: C_0(\omg^{-1}) \hookrightarrow LUC(\omg^{-1})$. 
Then $\jgt$ is bounded if and only if $L_{\gamma, \theta} < \infty$ and in this case, by taking dual maps, we obtain the commuting diagram: 

$$\xymatrixrowsep{3pc} \xymatrixcolsep{7pc}
\xymatrix{LUC(\omfi)^*  \ar@<.5ex>[r]^{\Jgt^*}  \ar@<.5ex>[rd]^{\jgt^*}  &   LUC(\omgi)^* \ar@<.5ex>[d]^{R_{C_0}=\iota_{C_0}^*} &
  \\
\M(\omf) \ar@{^{(}->}[u]^{\Theta_{LUC}} \ar@<.5ex>[r]^{\jgt^*\vert_{\M(\omf)}}    &
\M(\omg)     }$$
Observe that $\jgt^*: LUC(\omfi)^* \ra \M(\omg)$ is a weak$^*$-continuous algebra homomorphism. Abusing notation slightly, we let $\jgt^*$ also denote $\jgt^*\large{\vert}_{\M(\omf)} = \jgt^* \circ \Theta_{LUC}$. Obviously $\|\jgt^*\| \leq L_{\gamma, \theta}$ and $\sup_{x \in F} \|\jgt^*(\delta_x/\omf(x))\|_\omg = L_{\gamma, \infty}$, so  $\|\jgt^*\| = L_{\gamma, \theta}$.  When $\theta: F \ra G$ is proper, $\jgt$ maps $C_{00}(G)$ into $C_{00}(F)$ and therefore $\jgt$ maps $C_0(\omgi)$ into $C_0(\omfi)$ when $L_{\gamma, \theta}< \infty$. 

\bc  \label{jgt* a NM homom Corollary} Suppose that $L_{\gamma, \theta} < \infty$. Then $\jgt^*: \M(\omf) \ra \M(\omg)$ is a $\NM$-homomorphism with $\| \jgt^*\| = L_{\gamma,\theta}$, and $\jgt^*$ is the unique $so_l^\omf$-$\sigma^\omg$ continuous linear map $\M(\omf) \ra \M(\omg)$ mapping  $\delta_x$ to $ \gamma(x) \delta_{\theta(x)}$. When $\theta$ is proper, $\jgt^*: \M(\omf) \ra \M(\omg)$ is weak$^*$-continuous and when  $\gamma$ maps into $(0,\infty)$, $\jgt^*$ is positive.  \ec

\br  \label{M_gamma j_theta Remark}  \rm Let $1_F: F \ra \Cx$, ${\rm id}_F: F \ra F$ and let $M_\gamma := J_{\gamma, {\rm id}_F}$  $J_\theta := J_{1_F, \theta}$, $j_\theta := j_{1_F, \theta}$; this agrees with the notation used in the non-weighted case in \cite{Sto1}. From the general case, $M_\gamma: LUC(\omfi) \ra LUC (\omfi)$ is well-defined and bounded if and only if $ L_{\gamma, {\rm id}_F} = \sup_{x\in F} |\gamma(x)| < \infty$, and therefore if and only if $\gamma: F \ra \T$.  In this case, $M_\gamma$ is an isometric isomorphism of  $C_0(\omfi)$ onto $C_0(\omfi)$ with inverse $M_{\overline{\gamma}}$, and $M_\gamma^*: \M(\omf) \ra \M(\omf)$ is the unique weak$^*$-continuous isometric algebra  isomorphism mapping $\delta_x$ to $\gamma(x) \delta_x$ and, for a compact subgroup $K$ of $F$,  $m_K$ to $\gamma \vert_K m_K$.

  Also from the general case, $J_\theta$, $j_\theta$ are well-defined and bounded if and only if $L_\theta:= L_{1_F, \theta} = \sup_{x \in F} \omg(\theta(x))/\omf(x) < \infty$; in this case $\|J_\theta\| = \|j_\theta\| = L_\theta$, $J_\theta^*: LUC(\omfi)^* \ra LUC(\omgi)^*$ and $j_\theta^*: \M(\omf) \ra \M(\omg)$ are respectively the unique weak$^*$-continuous and $so_l^\omf-\sigma^\omg$ continuous algebra homomorphisms mapping $\delta_x$ to $\delta_{\theta(x)}$. 
  
When both $M_\gamma$ and $J_\theta / j_\theta$ are bounded,  $\Jgt = M_\gamma \circ J_\theta / \jgt = M_\gamma \circ j_\theta$. However, it is possible to have neither $M_\gamma$ nor $J_\theta/j_\theta$ bounded and yet $\Jgt/\jgt$ is bounded; see Remark \ref{M_gammax, j_thetaH unbounded Remark} below. 

\er 

Let $H$ be a closed subgroup of $G$ and define a continuous weight $\omx= \om_{\Cx \times H}$  on $\Cx \times H$ by putting $\omx(\alpha, s) = | \alpha|\omg(s)$. Let 
$$\gamma_\times = \gamma_{\Cx} : \Cx \times H \ra \Cx \quad \text{and} \quad \theta_H: \Cx \times H \ra H\leq G$$ be the projection homomorphisms. In the last statement of the following lemma, we further assume that $K$ is a compact  subgroup of $H$ with $\omg\equiv 1$ on $K$, $\rho \in \widehat{K}^1$ and 
$\Omega_\rho = \{ (\rho(k), k) : k \in K\}$. 

\blem \label{j_gammax,thetaH Lemma}   The map $j_{\gamma_\times, \theta_H}^* : \M(\om_\times) \ra \M(\omg)$ is a well-defined, weak$^*$-continuous contractive $\NM$-homomorphism  such that $j_{\gamma_\times, \theta_H}^*\delta_{(\alpha, s)} = \alpha \delta_s$. Moreover, $j_{\gamma_\times, \theta_H}^* (m_{\Omega_\rho}) = \rho m_K$. 
\elem 

\begin{proof}  For $(\alpha, s) \in \Cx \times H$, $|\gamma_\times(\alpha, s) |\omg(\theta_H(\alpha, s)) /  \omx(\alpha, s) =1$, so $\|j_{\gamma_\times, \theta_H} \| = L_{\gamma_\times, \theta_H} = 1$. Let $\psi \in C_0(\omg^{-1})$, $\ve >0$.  If $B$ is  a compact subset of $G$ such that $|\psi(s) / \omg(s)| < \ve $ for $s \in G \bs B$, then $|j_{\gamma_\times, \theta_H} \psi(\alpha, s) / \omx(\alpha, s)| = |\psi(s) / \omg(s)| < \ve $ for $(\alpha, s) \in  \Cx \times G \bs \{1\} \times B$; hence $j_{\gamma_\times, \theta_H}$ maps $C_0(\omg)$ into  $C_0(\omx)$. (Note though that $\theta_H$ is not proper.) By Corollary \ref{jgt* a NM homom Corollary}, only the last statement still requires justification. For this, note that $m_{\Omega_\rho} \in \M_{cr}(\Cx \times H) \subseteq \M_r(\Cx \times H) \cap \M(\omx)$, and since $\omg \equiv 1$ on $K$, $\omx \equiv 1$ on $\Omega_\rho$, and therefore  $\M(K, \omg) = \M_r(K) = M(K)$ and $\M(\Omega_\rho, \omx) =  \M_r(\Omega_\rho) = M(\Omega_\rho)$; the final statement hence follows from the argument used to establish \cite[Proposition 5.9]{Sto1}.
\end{proof}

\br \label{M_gammax, j_thetaH unbounded Remark} \rm   Although $\|j_{\gamma_\times, \theta_H}\|_{\infty, \omx^{-1}} =1$, neither $M_{\gamma_\times}$ nor $j_{\theta_H}$ is ever bounded because $\gamma_\times$ does not map into $\T$ and  $\om_G(\theta_H(\alpha, s)) / \omx(\alpha, s) = |\alpha|^{-1}$, so    $ L_{\theta_H} = \infty$.   \er

For the remainder of this subsection, $(H, \om)$ is a weighted locally compact group such that $\om\equiv 1$ on $K$, a compact normal subgroup of $H$, and  $\omhk (xK) := \om(x)$ is the continuous weight function on $H/K$ discussed in  Remark \ref{Constant on Cosets Remark Etc}.  Define $S_K: LUC(\omi) \ra \ell^\infty(\omhk^{-1})$ by putting 
$$S_K \psi(xK) = \int_K \psi(xk) \, d m_K(k) \qquad \text{for}  \ \psi \in LUC(\omi);$$
when we wish to emphasize that we are working specifically with $\om$ or the trivial weight, we will use the notation $S_K^\om$ or $S_K^1$, respectively. 

\bp \label{SK* Prop} 1. The linear map $S_K$  is a well-defined quotient map of $LUC(\omi)$ onto $LUC(\omhk^{-1})$ and of $C_0(\omi)$ onto $C_0(\omhk^{-1})$. 

\smallskip 

\noindent 2. The dual map $S_K^*$ is the unique weak$^*$-continuous isometric algebra isomorphic embedding of $LUC(\omhk^{-1})^*$ into $LUC(\omi)^*$ and of $\M(\omhk)$ into $\M(\om)$, mapping $\delta_{xK}$ to $\delta_x * m_K$. Moreover, $S_K^*: \M(\omhk) \hookrightarrow \M(\om)$ is a positive $\NM$-homomorphism. 
\ep

\begin{proof} 1. Let $\psi \in LUC(\omi)$. By \cite[Proposition 5.3 and its proof]{Sto1}, $S_K^1$ is a quotient mapping of $LUC(H)$ onto $LUC(H/K)$ (and of $C_0(H)$ onto $C_0(H/K)$), so $S_K^1(\psi/\om) \in LUC(H/K)$. Since 
$ S_K^1(\psi/\om) (xK) = (1 / \omhk(xK)) S_K^\om(\psi)(xK)$, $S_K^\om \psi \in LUC(\omhk^{-1})$ and 
$$ \| S_K^\om \psi\|_{\infty, \omhk^{-1}} = \left\| S_K^1\left(\psi / \om\right)\right\|_\infty \leq \left\|\psi / \om\right\|_\infty = \| \psi \|_{\infty, \omi}.$$
Taking $\phi$ in the unit ball of $LUC(\omhk^{-1})$, we can find $\psi_0 \in LUC(H)$ with $\| \psi_0\|_\infty \leq 1$ such that $S_K^1 \psi_0 = \phi/\omhk$; hence $\psi = \om \psi_0 \in LUC(H)$ with $\|\psi \|_{\infty, \omi} \leq 1$ and $S_K^\om \psi = \omhk S_K^1\psi_0 = \phi$. Thus, $S_K^\om$ is a quotient map of $LUC(\omi)$ onto $LUC(\omhki)$.  Observe that $\s(S_K^\om\psi) \subseteq q_K(\s(\psi) K^{-1})$ where $q_K:H \ra H/K$, so $S_K^\om$ maps $C_{00}(H)$ into $C_{00}(H/K)$, and therefore maps $C_0(\omi)$ into $C_0(\omhki)$. The above argument now shows that $S_K^\om$ is a quotient map of $C_0(\omi)$ onto $C_0(\omhki)$. 

\smallskip 

\noindent 2. We will establish the statement for $S_K^*: \M(\omhk) \hookrightarrow \M(\om)$; the same argument works for  $S_K^*: LUC(\omhki)^* \hookrightarrow LUC(\omi)^*$ (though we will not use this going forward). From part 1, $S_K^*$ is a weak$^*$-continuous linear isometric embedding. For $x \in H$ and $\psi \in C_0(\omi)$, 
$$\l S_K^*(\delta_{xK}), \psi\r_\om = \int_K\psi(xk) \, dm_K(k) = \iint \psi(st) \, d \delta_x(s) \, d m_K(t) = \l \delta_x * m_K, \psi \r_\om.$$
By \cite[Lemma 5.2]{Sto1}, $m_K$ is a central idempotent in $\M_{cr}(H)\subseteq \M(\om)$, so $S_K^*$ is a homomorphism by Lemma \ref{WK* continuous Homom Lemma}. Since $S_K^*(\Delta_{H/K}) \subseteq \M_{cr}(H)$ is a contractive subgroup of $\M_r(H)$, $S_K^*$ is a $\NM$-homomorphism, and $S_K^*$ is positive 
 because  $S_K$ maps $C_0(\omi)^+$ into $C_0(\omhki)^+$. 
\end{proof}  

Following the identification described in  Remark \ref{Remark mu_e}.2, when $H$ is a closed subgroup of $(G, \omg)$, but $K$ is not necessarily normal in $G$, we will sometimes abuse notation slightly and write $S_K^*: \M(\om_{H/K}) \hookrightarrow \M_r(\omg)$ in place of the homomorphism $R_H^* \circ S_K^*$. 

\subsection{The main theorem and some corollaries} 

\bt \label{Main NM Homom Thm} Let $\vp: \L^1(\omf) \ra \M(\omg)$ be a mapping. The following statements are equivalent: 
\bi \item[(i)] $\vp$ is a $\NM$-homomorphism and $\| \vp(\delta_{e_F})\|_{\omg} =1$; 
\item[(ii)] $\vp$ is a $\NM$-homomorphism and $\omg \equiv 1$ on $\s(\vp(\delta_{e_F})$;
\item[(iii)] there is a closed subgroup $H$ of $G$, a compact normal subgroup $\Omega_0$ of $\Cx \times H$ such that $\omx  \equiv 1$ on $\Omega_0$, and a continuous homomorphism $\theta:F \ra \Cx \times H / \Omega_0$ such that $\ds L_\theta = \sup_{x \in F} \omega_{\times/\Omega_0}(\theta(x))/  \omf(x) < \infty$ and $\vp = j_{\gamma^\times, \theta_H}^* \circ S_{\Omega_0}^* \circ j_\theta^*$ i.e., $\vp$ factors as \beq  \label{Main Thm Commuting Diagram Eqn}  \xymatrixrowsep{2pc} \xymatrixcolsep{4pc}
\xymatrix{\L^1(\omf) \ar@<.1ex>[r]^{j_\theta^* \quad  } & \M(\omega_{\times/\Omega_0})    \ar@{^{(}->}[r]^{S_{\Omega_0}^*} &\M(\omx)  \ar@<.1ex>[r]^{ j_{\gamma^\times, \theta_H}^*}  &   \M(\omg).  } \eeq

\ei
In statement (iii), $\| \vp \| = \|j_\theta^*\| = L_\theta$ and we can take $H = \s(\Gamma_\vp)$, $\vp(\delta_{e_F}) = \rho m_K$ and $\Omega_0 = \Omega_\rho$ associated to $\Gamma_\vp$ as in Lemma \ref{Main NM-subgroup Description Lemma}.
\et

\begin{proof} (i) $\Rightarrow$ (ii) This is an immediate consequence of  Corollary \ref{Norm-one idempotents in M(w) Corollary}.

\smallskip 

\noindent (ii) $\Rightarrow$ (iii)  Since $\Gamma = \Gamma_\vp$ is a $\NM$-subgroup of $\M_r(G)$, employing Proposition \ref{NM-subgroup Isomorphism Prop} and the notation therein,  $$\theta: F \ra \Cx \times H/\Omega_\rho: x \mapsto \phi^{-1}(\vp(\delta_x))$$ is a continuous homomorphism, where 
$$\phi: \Cx \times H / \Omega_\rho \ra \Gamma_{\Cx \times H } \supseteq \Gamma: (\alpha, s)\Omega_\rho \mapsto \alpha \delta_s * \rho m_K.$$
Note that $\omx \equiv 1$ on $\Omega_\rho$ since $\omg\equiv 1$ on $K = \s(\iota_\Gamma)$, and therefore $\om_{\times/\Omega_\rho}$ is a well-defined weight on $\Cx \times H/\Omega_\rho$.   Let $x \in F$ and suppose that $\vp(\delta_x) = \alpha \delta_s * \rho m_K$. Then $\theta(x) = (\alpha, s) \Omega_\rho$ and $\|\vp(\delta_x) \|_\omg = \omg(s) |\alpha|$ by Remark \ref{Properties of alpha deltas * rho mK Remark}.2, so 
$${\om_{\times/\Omega_\rho}(\theta(x))\over \omf(x)} = {\omx(\alpha, s) \over \omf(x) } = {|\alpha|\omg(s)\over \omf(x)}  = \left\|  \vp\left({\delta_x \over \omf(x)}\right)\right\|_\omg \leq \|\vp \| \left\|  {\delta_x \over \omf(x)} \right\|_\omf = \| \vp \|.$$
Hence, $L_\theta = \sup_{x \in F} \om_{\times/\Omega_\rho}(\theta(x))/\omf(x) \leq \| \vp \|$, so $j_\theta^*$ is a bounded (with $\|j_\theta^*\|\leq \|\vp\|$)  $so_l^\omf-\sigma^{\om_{\times/\Omega_\rho}}$ continuous homomorphism, and by Lemma \ref{j_gammax,thetaH Lemma} and Proposition \ref{SK* Prop}, $j_{\gamma_\times, \theta_H}^* \circ S_{\Omega_\rho}^* \circ j_\theta^*$ is $so_l^\omf-\sigma^\omg$ continuous. Supposing again that $\vp(\delta_x) = \alpha \delta_s * \rho m_K$, 
\beq \label{Main Thm Calculation Eqn} \begin{aligned} j_{\gamma_\times, \theta_H}^* \circ S_{\Omega_\rho}^* \circ j_\theta^* (\delta_x) & =  j_{\gamma_\times, \theta_H}^* \circ S_{\Omega_\rho}^* (\delta_{(\alpha, s) \Omega_\rho} ) =  j_{\gamma_\times, \theta_H}^* (\delta_{(\alpha, s)}*m_{\Omega_\rho}) \\ &=  j_{\gamma_\times, \theta_H}^* (\delta_{(\alpha, s)})* j_{\gamma_\times, \theta_H}^*(m_{\Omega_\rho}) = \alpha \delta_s * \rho m_K  = \vp (\delta_x), \end{aligned} \eeq  
where we have used Lemma \ref{j_gammax,thetaH Lemma} once more.  By $so_l^\omf$-density of $\C F$ in $\M(\omf)$, $\vp =  j_{\gamma_\times, \theta_H}^* \circ S_{\Omega_\rho}^* \circ j_\theta^*$; since  $j_{\gamma_\times, \theta_H}^*$ and  $S_{\Omega_\rho}^*$ are contractions, $\|\vp\| =  \| j_\theta^* \| = L_\theta$.   

\smallskip 

\noindent (iii) $\Rightarrow$ (i) By Corollary \ref{jgt* a NM homom Corollary}, Lemma \ref{j_gammax,thetaH Lemma} and Proposition \ref{SK* Prop}, $\vp = j_{\gamma_\times, \theta_H}^* \circ S_{\Omega_0}^* \circ j_\theta^*$ is a $so_l^\omf-\sigma^\omg$ continuous composition of $\NM$-homomorphisms. By Proposition \ref{Composition of NM-homs Prop}, $S_{\Omega_0}^* \circ j_\theta^*$ is a $\NM$-homomorphism because $ S_{\Omega_0}^* \circ j_\theta^*(\delta_{e_F}) = m_{\Omega_0} \in \M_r(\Cx \times H)$ and $\|m_{\Omega_0}\| = 1$. Since $\omx \equiv 1$ on $\Omega_0$, $\|m_{\Omega_0}\|_\omx = \| m_{\Omega_0}\| =1$, and therefore $\| \vp(\delta_{e_F})\|_\omg = \| j_{\gamma_\times, \theta_H}^*(m_{\Omega_0})\|_\omg \leq 1$ because $j_{\gamma_\times, \theta_H}^*$ is a contraction. By Corollary \ref{Norm-one idempotents in M(w) Corollary}, $ \vp(\delta_{e_F})=   j_{\gamma_\times, \theta_H}^* \circ (S_{\Omega_0}^* \circ j_\theta^*)(\delta_{e_F}) \in \M_r(G)$ with $\| \vp(\delta_{e_F})\| =1$, so, by a second application of Proposition \ref{Composition of NM-homs Prop}, $\vp$ is a $\NM$-homomorphism.
\end{proof}  

\bex \label{NM vs NMOG Homom} \rm There are $\NM$-homomorphisms $\vp: \L^1(\omf) \ra \M(\omg)$ for which $\| \vp(\delta_{e_F})\|_\omg \neq 1$. To see this, take $K$ to be any example of  compact subgroup of $G$ on which $\omg$ is not identically $1$. Then $\iota = m_K \in \M_r(G) \cap \M(\omg)$ is an idempotent with $\|\iota\| =1$, but $\| \iota \|_\omg \neq 1$ by Corollary \ref{Norm-one idempotents in M(w) Corollary}. Let $F$ be a finite group and define a homomorphism  $\vp: \L^1(F) = \C F \ra \M(\omg)$ by putting $\vp(f) = \sum_{s\in F} f(s) \, \iota$. Then $\Gamma_\vp = \vp(\Delta_F) = \{\iota\}$, so $\vp$ is a $\NM$-homomorphism, but $\|\vp(\delta_{e_F})\|_\omg \neq 1$ and (therefore) $\vp$ is not a $\NM_\omg$-homomorphism.   
\eex

\bc \label{NMOG-Homom Corollary} Let $\vp: \L^1(\omf) \ra \M(\omg)$ be a mapping. The following statements are equivalent: 
\bi \item[(i)] $\vp$ is a $\NM_\omg$-homomorphism; 
\item[(ii)] there is a closed subgroup $H$ of $G$ on which $\omg$ is multiplicative, a compact normal subgroup $\Omega_0$ of $\Cx \times H$, and a continuous homomorphism $\theta:F \ra \Cx \times H / \Omega_0$ such that $ \ds L_\theta = \sup_{x \in F} \omega_{\times/\Omega_0}(\theta(x))/  \omf(x) < \infty$ and $\vp = j_{\gamma^\times, \theta_H}^* \circ S_{\Omega_0}^* \circ j_\theta^*$ i.e., $\vp$ has the factorization from (\ref{Main Thm Commuting Diagram Eqn}). 
\ei 
In statement (ii), $\| \vp \| = \|j_\theta^*\| = L_\theta$ and we can take $H = \s(\Gamma_\vp)$, $\vp(\delta_{e_F}) = \rho m_K$ and $\Omega_0 = \Omega_\rho$ associated to $\Gamma_\vp$ as in Lemma \ref{Main NM-subgroup Description Lemma}.
\ec

\begin{proof} Assume statement  (i) holds. As noted in Remark \ref{NM vs NM_omega Homoms Remark}.2(i), $\vp$ is a $\NM$-homomorphism with $\| \vp(\delta_{e_F})\|_\omg =1$, so statement (iii) of Theorem \ref{Main NM Homom Thm} holds, with  $H = \s(\Gamma_\vp)$. By Theorem \ref{NMO-subgroup Description Theorem}, $\omg$ is multiplicative on $H$. Conversely, assume statement (ii) holds. Then  $\omx$ is multiplicative on $\Cx \times H$, so $\omx \equiv 1$ on the compact subgroup $\Omega_0$ of $\Cx \times H$. By Theorem \ref{Main NM Homom Thm}, $\Gamma_\vp = \vp(\Delta_F)$ is a $\NM$-subgroup of $\M_r(G)$. Observe that the weak$^*$-continuous homomorphism $j_{\gamma_\times, \theta_H}^*$ maps $\C(\Cx\times H)$ into the weak$^*$-closed subalgebra $\M_H(\omg)=\{ \mu \in \M(\omg): \s(\mu) \subseteq H\}$ --- see Lemma \ref{j_gammax,thetaH Lemma} and Remark \ref{Remark mu_e}.2 --- and therefore maps $\M(\omx)$ into $\M_H(\omg)$. Hence, $\omg$ is multiplicative on $\s(\Gamma_\vp) \subseteq H$ and by Theorem \ref{NMO-subgroup Description Theorem} we conclude that $\vp$ is a $\NM_\omg$-homomorphism.  
\end{proof} 

\br \rm When $\omf \equiv 1$ and $\vp: \L^1(F) \ra \M(\omg)$ is a contractive homomorphism, $\vp$ is a $\NM_{\omg}$-homomorphism, so Corollary \ref{NMOG-Homom Corollary} extends \cite[Theorem 5.11]{Sto1}, which describes all contractive homomorphisms $\vp:\L^1(F) \ra \M_r(G)$ --- the case when $\omg$ is also trivial.  Explicitly, see below. (Note, however, that: we have made significant use herein of several arguments from \cite{Sto1}; and  for non-trivial weights, even $\NM$-{\it isomorphisms} can fail to be contractive --- see  \cite[Example 3.9]{Gha-Zad}.)
\er 

\bc \label{Sto1 Thm Corollary} Let $\vp:\L^1(F) \ra \M_r(G)$ be a map. The following statements are equivalent: 
\bi \item[(i)]  $\vp$ is a contractive homomorphism; 
\item[(ii)] $\vp$ is a $\NM$-homomorphism; 
\item[(iii)]  $\vp$ has the factorization described in \cite[Theorem 5.11]{Sto1}.
\ei 
\ec 

\begin{proof} Only (ii) implies (iii) requires justification, for which we employ Corollary \ref{NMOG-Homom Corollary}: Since $\omx(\alpha, s)= |\alpha|$ is a continuous homomorphism containing $\Omega_\rho$ in its kernel, $\om_{\times/\Omega_\rho}$ is a continuous homomorphism in this case. As $L_\theta < \infty$, $\om_{\times/\Omega_\rho} \circ \theta$ is a bounded, and therefore trivial, homomorphism into $(0, \infty)$; it follows that $\theta$ maps into $\T \times H / \Omega_\rho$, on which $\om_{\times/\Omega_\rho} \equiv 1$. Thus, Corollary \ref{NMOG-Homom Corollary} yields  the desired factorization (with trivial weights, $\gamma_\times: \T \times H \ra \T$, and $j_{\gamma_\times, \theta_H}^* = j_{\theta_H}^* \circ M_{\gamma_\times}^*$ by Remark \ref{M_gamma j_theta Remark}).
\end{proof}

The following terminology is inspired by the   notion of standard isomorphism in \cite{Gha-Zad}.


\bd \rm A pair $(\gamma, \theta)$ will be called a \it standard homomorphism of $(F, \omf)$ into $(G, \omg)$ \rm if $$\gamma: F \ra \Cx \qquad \text{and} \qquad \theta: F \ra G$$ are continuous homomorphisms  such that  $\ds L_{\gamma, \theta}< \infty$; 
 when $\gamma$ maps into the positive reals,  $(0, \infty)$, we will call $(\gamma,\theta)$ a \it positive standard homomorphism.\rm 
\ed

The following result and its proof are weighted analogues of \cite[Corollary 5.8]{Sto1} and its proof. 

\bc \label{When rho extends to H Corollary}  Let $\vp: \L^1(\omf) \ra \M(\omg)$ be a $\NM$-homomorphism such that $\| \vp(\delta_{e_F})\|_\omg =1$, written as $\vp(\delta_{e_F}) = \rho m_K$ with $\omg\equiv 1$ on $K$ as in Corollary \ref{Norm-one idempotents in M(w) Corollary}. If $\rho$ extends to $\rho_H \in \widehat{H}^1$ --- e.g., when $H$ is abelian --- where $H = \s(\Gamma_\vp)$, then there is a standard homomorphism $(\gamma, \theta_K)$  of $(F, \omf)$ into $(H/K, \om_{H/K})$ such that  $\vp= M_{\rho_H}^* \circ S_K^* \circ j_{\gamma, \theta_K}^*$. That is, $\vp$ factors as  
\beqs    \xymatrixrowsep{2pc} \xymatrixcolsep{4pc}
\xymatrix{ \L^1(\omf)\subseteq \M(\omf) \ar@<.5ex>[r]^{j_{\gamma, \theta_K}^*} & \M(\om_{H/K})   \ar@{^{(}->}[r]^{S_K^*} & \M(\om_H)  \ar@{^{(}->}[r]^{M_{\rho_H}^*}  &   \M(\omg).  } \eeqs
\ec 

\begin{proof}    By Corollary \ref{Norm-one idempotents in M(w) Corollary}, $\om_K \equiv 1$  on $K$, so $\om_{H/K}$ is well-defined and by Theorem \ref{Main NM Homom Thm} $\vp$ has the factorization (\ref{Main Thm Commuting Diagram Eqn}), where $\theta: F \ra \Cx \times H/\Omega_\rho$. One can check that 
$$p_{\Cx}: \Cx \times H / \Omega_\rho \ra \Cx: (\alpha, s) \Omega_\rho \mapsto \alpha \overline{\rho_H(s)} \quad \text{and} \quad p_K: \Cx \times H / \Omega_\rho \ra H/K: (\alpha, s) \Omega_\rho \mapsto sK$$ are well-defined continuous homomorphisms, so $\gamma:= p_{\Cx} \circ \theta$ and $\theta_K := p_K \circ \theta$ are continuous homomorphisms. If $\theta(x) = (\alpha, s) \Omega_\rho$, then 
$$|\gamma(x)|\om_{H/K}(\theta_K(x)) = |\alpha| | \overline{\rho_H(s)}|\om_{H/K}(sK) = |\alpha|\om(s) = \om_{\times/\Omega_\rho}(\theta(x))$$
so 
$$\|j_{\gamma, \theta_K}\| = L_{\gamma, \theta_K} = \sup_{x \in F} |\gamma(x)| {\om_{H/K}(\theta_K(x))\over \omf(x)}  = \sup_{x\in F} {\om_{\times/\Omega_\rho}(\theta(x)) \over \omf(x)}  = L_\theta = \| \vp\|.$$
When $\theta(x) = (\alpha, s) \Omega_\rho$,  the calculation (\ref{Main Thm Calculation Eqn}) shows that $\vp(\delta_x) =  j_{\gamma_\times, \theta_H}^* \circ S_{\Omega_\rho}^* \circ j_\theta^* (\delta_x) = \alpha \delta_s * \rho m_K$ and 
\beqs M_{\rho_H}^* \circ S_K^* \circ j_{\gamma, \theta_K}^*(\delta_x) & = &  = M_{\rho_H}^* \circ S_K^* (\alpha \overline{\rho_H(s)} \delta_{sK})    =  \alpha \overline{\rho_H(s)} M_{\rho_H}^*(\delta_s * m_K) \\   
 & = &  \alpha \overline{\rho_H(s)}  \rho_H(s) \delta_s * (\rho_H\vert_K)m_K 
  =  \alpha \delta_s * \rho m_K = \vp (\delta_x); 
\eeqs 
by $so_l^\omf-\sigma^\omg$ continuity of the maps involved, the proof is complete. 
\end{proof}

\bc \label{phi(delta_e) positive Corollary} Let $\vp: \L^1(\omf) \ra \M(\omg)$ be a mapping. Then:
\bi \item[(a)] $\vp$ is a $\NM$-homomorphism such that $\vp(\delta_{e_F})\geq 0$ and  $\| \vp(\delta_{e_F})\|_\omg =1$  if and only if  
there is  a standard homomorphism $(\gamma, \theta_K)$ of $(F, \omf)$ into $(H/K, \om_{H/K})$ such that $\vp= S_K^* \circ j_{\gamma, \theta_K}^*$, where $H$ is a closed subgroup of $G$ containing the compact normal subgroup $K$ on which $\omg\equiv 1$;  moreover we can take $H=\s(\Gamma_\vp)$. 
\item[(b)] $\vp$ is a $\NM$-homomorphism such that $\vp(\delta_{e_F}) = \delta_{e_G}$   if and only if there is a standard homomorphism $(\gamma, \theta)$  of $(F, \omf)$ into $(G, \omg)$  such that  $\vp=  j_{\gamma, \theta}^*$.\ei 
\ec

\begin{proof}  For the forward implication in (a), $\vp(\delta_{e_F}) = m_K$ for some compact subgroup $K$  of $G$  by Theorem \ref{Positive Idempotent Thm}, so we can take  $\rho_H \equiv 1$ in Corollary \ref{When rho extends to H Corollary}. The converse direction in (a) follows from Proposition \ref{Composition of NM-homs Prop}  because $S_K^* \circ j_{\gamma, \theta_K}^*(\delta_{e_F}) = m_K$,  a positive norm-one idempotent in both $\M_r(G)$ and $\M(\omg)$ since $\omg \equiv 1$ on $K$.  Part  (b) follows from (a) because  $K=\{e_G\}$ in this case.
\end{proof} 

\bc \label{Positive Beurling algebra Homoms Corollary} Let $\vp: \L^1(\omf) \ra \M(\omg)$ be a mapping. Then:
\bi \item[(a)] $\vp$ is  a positive homomorphism such that $\| \vp(\delta_{e_F})\|_\omg =1$  if and only if  
 there is  a positive  standard homomorphism $(\gamma, \theta_K)$ of $(F, \omf)$ into $(H/K, \om_{H/K})$ such that $\vp= S_K^* \circ j_{\gamma, \theta_K}^*$, where $H$ is a closed subgroup of $G$ containing the compact normal subgroup $K$ on which $\omg\equiv 1$;  moreover we can take $H=\s(\Gamma_\vp)$.  
\item[(b)] $\vp$ is a positive homomorphism such that $\vp(\delta_{e_F}) = \delta_{e_G}$   if and only if  there is a positive standard homomorphism $(\gamma, \theta)$  of $(F, \omf)$ into $(G, \omg)$  such that  $\vp=  j_{\gamma, \theta}^*$.\ei 
\ec

\begin{proof} Only the forward implication of part (a) requires explanation. As noted in Remark \ref{Canonical extension to M(omf) Remark}, $\vp$ is a positive  homomorphism on $\M(\omf)$. By Proposition \ref{NMw-subgroup Prop}(iv), Corollary \ref{phi(delta_e) positive Corollary} provides everything except that $\gamma$ maps into the positive reals. However, if $\theta_K(x) = sK$, positivity of  $\vp(\delta_x) = S_K^* \circ j_{\gamma, \theta_K}(\delta_x) = \gamma(x) \delta_s * m_K$ implies that $\gamma(x) >0$.
\end{proof} 

The following important special case, which characterizes \it all \rm positive homomorphisms $\vp: \L^1(F) \ra \M_r(G)$, is new and worth explicitly recording. Bipositive isomorphisms between group algebras were first described by  Kawada in 1948 \cite{Kaw}.

\bc \label{Positive homs unweighted Corollary}  (a) A map $\vp: \L^1(\omf) \ra \M_r(G)$ is a  positive homomorphism  if and only if  
there is a compact normal subgroup $K$ of a closed subgroup $H$  of $G$ and continuous homomorphisms $\gamma: F \ra (0,\infty)$  and $\theta_K : F \ra H/K$ such that  $\sup_{x \in F} \gamma(x)/ \omf(x) < \infty$ and $\vp= S_K^* \circ j_{\gamma, \theta_K}^*$.

\smallskip 

\noindent  (b)  A map $\vp: \L^1(F) \ra \M_r(G)$ is a  positive homomorphism  if and only if  
there is a compact normal subgroup $K$ of a closed subgroup $H$  of $G$ and a continuous homomorphism $\theta_K : F \ra H/K$ such that  $\vp= S_K^* \circ j_{\theta_K}^*$. Thus, bounded positive homomorphisms $\L^1(F) \ra \M_r(G)$  are automatically contractive. 
\ec 

\begin{proof} Let $\vp: \L^1(\omf) \ra \M_r(G)$ be a positive homomorphism. By Prop. \ref{NMw-subgroup Prop} (iv) and (iii), $\Gamma_\vp$ is a $\NM$-subgroup of $\M_r(G)$ and $\|\vp(\delta_{e_F})\| =1$. Hence, Corollary \ref{Positive Beurling algebra Homoms Corollary} applies with $\ds L_{\gamma, \theta_K} = \sup_{x\in F} \gamma(x) / \omf(x) < \infty$ and $\vp =  S_K^* \circ j_{\gamma, \theta_K}^*$. When $\omf$ is trivial, $\ds  \sup_{x \in F} \gamma(x) < \infty$ implies  that  $\gamma \equiv 1$. 
\end{proof} 

\subsection{Recovery of the isomorphism theorems} 

We now recover the main isomorphism theorems from \cite{Gha-Zad, Kaw, Zad}, and improve the results from \cite{Gha-Zad} and \cite{Kaw} by showing that positive isomorphisms are automatically bipositive. Our proofs are almost entirely independent of \cite{Gha-Zad} and \cite{Zad}, however we will use the following lemma, in which part (a) is Zadeh's important observation \cite[Lemma 2.2]{Zad}  and part (b) is \cite[Remark 2.14]{Sto2}, which  employed \cite[Lemma 3.3]{Gha-Zad}; the proofs of these lemmas   can be read independently of the rest of \cite{Zad} and \cite{Gha-Zad}. In the statements that follow, isomorphisms are always assumed to be surjective.  

\blem \label{Zadeh Ghar-Zadeh Lemma} (i) The set of extreme points of the unit ball of  $\M(\omg)$ is $\left\{ {\alpha\over \omg(x)} \, \delta_x : \alpha \in \T, \  x \in G\right\}$.

\smallskip 
\noindent (ii) If $\vp: \M(\omf) \ra \M(\omg)$ is a Banach algebra isomorphism, then it is $so_l^\omf - \sigma^\omg$ continuous; it is thus the unique $so_l^\omf - \sigma^\omg$-continuous extension of $\vp\vert_{\L^1(F)}$ to $\M(\omf)$.  
\elem 

For clarity,  in the next lemma we use $\ov{\vp}$ to denote the $so_l^\omf - \sigma^\omg$ continuous extension to $\M(\omf)$  of  a homomorphism $\vp$ on $\L^1(\omf)$. 

\blem \label{Extension of isom from L1(w) to M(w) Lemma} Let $\vp: \L^1(\omf) \ra \L^1(\omg)$ be a Banach algebra isomorphism. Then $\ov{\vp}: \M(\omf) \ra \M(\omg)$ is a $so_l^\omf-so_l^\omg$ continuous Banach algebra isomorphism and $\ov{\vp}^{-1} = \ov{\vp^{-1}}$. Moreover, $\ov{\vp}$ is an isometry when $\vp$ is an isometry.
\elem

\begin{proof}  Suppose that $\mu_i \ra \mu$ $so_l^\omf$ and take $g \in \L^1(\omg)$. Taking $f \in \L^1(\omf)$ such that $\vp(f) = g$, $\| g* (\ov{\vp}(\mu_i) - \ov{\vp}(\mu)\|_\omg \leq \|\vp\|\|f * \mu_i - f * \mu\|_\omf \ra 0,$ giving $so_l^\omf-so_l^\omg$ continuity of $\ov{\vp}$. (Note that only surjectivity of $\vp$ was required here.) Thus, $\ov{\vp^{-1}} \circ \ov{\vp}$ is $so_l^\omf - \sigma^\omf$ continuous on $\M(\omf)$, and agrees with $\text{id}_{\M(\omf)}$ on $\L^1(\omf)$, which is also $so_l^\omf - \sigma^\omf$ continuous, (e.g., by \cite[Proposition 2.9]{Sto2} with $\S(\omf^{-1}) = C_0(\omf^{-1})$); hence $\ov{\vp^{-1}} \circ \ov{\vp} = \text{id}_{\M(\omf)}$, and  $\ov{\vp} \circ \ov{\vp^{-1}}  = \text{id}_{\M(\omg)}$. If $\vp$ (and $\vp^{-1}$) is an isometry, then $\ov{\vp}$ and $\ov{\vp^{-1}} = \ov{\vp}^{-1}$ are both contractions by Remark \ref{Canonical extension to M(omf) Remark}, so  $\ov{\vp}$ is an isometry. 
\end{proof}

\bd \rm A pair $(\gamma, \theta)$ is called a \it standard isomorphism of $(F, \omf)$ onto $(G, \omg)$ \rm if $\gamma: F \ra \Cx$ is a continuous homomorphism and $\theta: F \ra G$ is a topological group isomorphism such that 
$$0< \inf_{s \in F} |\gamma(s)| {\omg(\theta(s))\over \omf(s)}=: l_{\gamma, \theta} \quad {\rm and } \quad L_{\gamma, \theta} = \sup_{s \in F} |\gamma(s)| {\omg(\theta(s))\over \omf(s)}< \infty;$$
when  $\gamma: F \ra (0, \infty)$, we call $(\gamma, \theta)$ a \it positive standard isomorphism;  \rm when $|\gamma| = \omf/\omg\circ \theta$ --- in which case $l_{\gamma,\theta} = L_{\gamma,\theta}=1$ --- we call $(\gamma, \theta)$ a \it linked standard isomorphism.\rm 
\ed 

\br \rm   What we have called a positive standard isomorphism is called a standard isomorphism in \cite{Gha-Zad}; in the terminology of \cite{Zad},  $\theta$ is an isomorphism of the weighted locally compact groups $(F, \omf)$ and $(G, \omg)$ when $(\gamma, \theta)$ is a linked isomorphism.  \er

\bp  \label{Easy Direction M(w) Prop}   Suppose that  $(\gamma, \theta)$ is a standard isomorphism of $(F, \omf)$ onto $(G, \omg)$, and let $\beta(t):=1/ \gamma(\theta^{-1}(t))$ for $t \in G$.  Then  $\jgt^*$ is a $\sigma^\omf - \sigma^\omg$ continuous Banach algebra isomorphism with $\vp^{-1} = j_{\beta, \theta^{-1}}^*$. Moreover,  $\jgt^*$ is a bipositive isomorphism when $(\gamma, \theta)$ is a positive standard isomorphism and $\jgt^*$ is an isometric isomorphism when $(\gamma, \theta)$ is a linked standard isomorphism. \ep 

\begin{proof} Observe that $L_{\beta,\theta^{-1}} \leq 1/l_{\gamma, \theta}< \infty$, so by Corollary \ref{jgt* a NM homom Corollary}, $\vp= \jgt^*$ and $j_{\beta, \theta^{-1}}^*$ are $\NM$-homomorphisms that are weak$^*$-continuous because $\theta$ and $\theta^{-1}$ are proper maps. For $s \in F$ and $t \in G$, $j_{\beta, \theta^{-1}}^* \circ \jgt^*(\delta_s) = \delta_s$ and 
$\jgt^* \circ j_{\beta, \theta^{-1}}^* = \delta_t $, so $\vp^{-1} = j_{\beta, \theta^{-1}}^*$ by Lemma \ref{WK* continuous Homom Lemma}. When $\gamma$ maps into $(0,\infty)$, so does $\beta$, and therefore both $\vp$ and $\vp^{-1}$ are positive. When $(\gamma, \theta)$ is a linked standard isomorphism, $L_{\gamma,\theta} = l_{\gamma, \theta} =1$, so both $\vp$ and $\vp^{-1}$ are contractions, and therefore $\vp$ is isometric. 
\end{proof} 

By Lemma \ref{Zadeh Ghar-Zadeh Lemma} (ii), when $\vp:\M(\omf) \ra \M(\omg)$ is a Banach algebra isomorphism the $so_l^\omf-\sigma^\omg$ continuity requirement  for  $\vp$ to be a $\NM$-homomorphism (Definition \ref{NM-Homoms Defn}) is automatically satisfied; this condition can thus be dropped from the definition of a  ``$\NM$-isomorphism" $\vp: \M(\omf) \ra \M(\omg)$.

\bt \label{biNM Isomorphism Theorem} Suppose that $\vp: \L^1(\omf) \ra \L^1(\omg)$ or $\vp: \M(\omf) \ra \M(\omg)$ is a $\NM$-isomorphism such that $\vp^{-1}$ is also a $\NM$-isomorphism.  Then $\vp = \jgt^*$ for some standard isomorphism $(\gamma, \theta)$
 of $(F, \omf)$ onto $(G, \omg)$. 
\et

\begin{proof}  By Lemma \ref{Extension of isom from L1(w) to M(w) Lemma}, it suffices to establish the theorem when $\vp:\M(\omf) \ra \M(\omg)$. In this case, by Lemma \ref{Zadeh Ghar-Zadeh Lemma} (ii) and Corollary \ref{phi(delta_e) positive Corollary} (b), there are standard homomorphisms $(\gamma, \theta)$ and $(\beta, \psi)$ of $(F, \omf)$ into $(G, \omg)$ and $(G, \omg)$ into $(F, \omf)$, respectively, such that $\vp = \jgt^*$ and $\vp^{-1} = j_{\beta,\psi}^*$. Hence, for  $s \in F$ and $t \in G$, $$\delta_s  = \vp^{-1}(\vp(\delta_s)) = \gamma(s) \beta(\theta(s)) \delta_{\psi(\theta(s))} \qquad \text{and} \qquad \delta_t  = \vp(\vp^{-1}(\delta_t)) = \beta(t) \gamma(\psi(t)) \delta_{\theta(\psi(t))}.$$
Thus, $\theta: F \ra G$ is a topological isomorphism with $\psi = \theta^{-1}$, $\gamma(s) = 1/\beta(\theta(s))$ and $\beta(t) = 1/\gamma(\theta^{-1}(t))$; from this it is easy to check that $l_{\gamma, \theta} \geq 1/L_{\beta, \theta^{-1}}>0$.  
\end{proof} 

 
 We now establish  Ghahramani and Zadeh's Theorems 3.5 and 4.3 of \cite{Gha-Zad}, the main results on Beurling group and measure algebras therein, with the weaker hypothesis of positivity replacing bipositivity. (If one assumes below that $\vp$ is bipositive, the result follows immediately from Theorem \ref{biNM Isomorphism Theorem} and Remark \ref{NM vs NM_omega Homoms Remark}.2 (iii).) 

\bt \label{Positive Isom Thm} Let $\vp: \L^1(\omf) \ra \L^1(\omg)$ or $\vp:\M(\omf) \ra \M(\omg)$ be a positive isomorphism. Then $\vp$ is bipositive and $\vp= \jgt^*$ for some positive standard isomorphism $(\gamma, \theta)$ of $(F, \omf)$ onto $(G, \omg)$. 
\et 

\begin{proof} By Remark  \ref{Canonical extension to M(omf) Remark}, $\vp$ is positive on $\M(\omf)$ when it is positive on $\L^1(\omf)$, so by Lemma \ref{Extension of isom from L1(w) to M(w) Lemma} it suffices to establish the theorem when $\vp: \M(\omf) \ra \M(\omg)$ is a positive isomorphism.  By Corollary \ref{Positive Beurling algebra Homoms Corollary} (b), there is a positive standard homomorphism $(\gamma, \theta)$ of $(F, \omf)$ into $(G, \omg)$ such that $\vp = \jgt^*$. For $s \in F$, $\vp(\delta_s) = \gamma(s) \delta_{\theta(s)}$, so $H = \s(\Gamma_\vp) = \ov{\theta(F)}$ and therefore $\vp$ maps $\C F$ into $\M_H(\omg)$, a $\sigma^\omg$-closed subalgebra of $\M(\omg)$ --- see Remark \ref{Remark mu_e}.2. Since $\C F$ is $so_l^\omf$-dense in $\M(\omf)$ and $\vp$ is $so_l^\omf-\sigma^\omg$ continuous by Lemma \ref{Zadeh Ghar-Zadeh Lemma} (ii), $\vp$ maps $\M(\omf)$ into $\M_H(\omg)$; hence, $G = H = \ov{\theta(F)}$. 

Let $t \in G$ and take a net $(s_i)$ in $F$ such that $\theta(s_i) \ra t$. Then $\delta_{\theta(s_i)} \ra \delta_t$ $so_l^\omg$, so by Lemma \ref{Zadeh Ghar-Zadeh Lemma}(ii), $\vp^{-1}(\delta_{\theta(s_i)}) \ra \vp^{-1}(\delta_t)$ $\sigma^\omf$. But for each $i$, $\vp^{-1}(\delta_{\theta(s_i)}) = \vp^{-1}(1/\gamma(s_i) \vp(\delta_{s_i}))= 1/\gamma(s_i)\delta_{s_i} \geq 0$, so $\vp^{-1}(\delta_t) \geq 0$ as well. Hence, $\Gamma_{\vp^{-1}} = \vp^{-1}(\Delta_G) \subseteq \M(\omf)^+$ and $\vp^{-1}(\delta_{e_G}) = \delta_{e_F}$, so by Remark \ref{NM vs NM_omega Homoms Remark}.2 (iii), $\vp^{-1}$ is a $\NM$-isomorphism. By Theorem \ref{biNM Isomorphism Theorem}, $(\gamma, \theta)$ is a  standard isomorphism of $(F, \omf)$ onto $(G, \omg)$, that we already know is positive. By Proposition \ref{Easy Direction M(w) Prop}, $\vp$ is bipositive.
\end{proof}

The following result was proved by S. Zadeh as Theorems 2.4, 2.5 and 3.4 of  \cite{Zad}. It should be noted that her use of the notation $\jgt^*$ is slightly different from ours.   We use her Lemma \ref{Zadeh Ghar-Zadeh Lemma} (i) and her observation that linear isometric isomorphisms preserve extreme points of the unit ball.

\bt \rm (Zadeh) \label{Isometric Isom Thm} Let $\vp: \L^1(\omf) \ra \L^1(\omg)$ or $\vp: \M(\omf) \ra \M(\omg)$ be an isometric isomorphism. Then $\vp = \jgt^*$ for some linked standard isomorphism $(\gamma, \theta)$ of $(F, \omf)$ onto $(G, \omg)$. 
\et 

\begin{proof} By Lemma 4.22, we only need to establish the corollary when $\vp: \M(\omf) \ra \M(\omg)$ is an isometric isomorphism. As a linear isometry, $\vp$ maps extreme points of the unit ball of $\M(\omf)$ to extreme points of the unit ball of $\M(\omg)$. Thus, by Lemma \ref{Zadeh Ghar-Zadeh Lemma} (i), for each $s \in F$, 
$\vp(\delta_s) = \omf(s) \vp (1/\omf(s) \, \delta_s) \in \{\alpha \delta_t : \alpha \in \Cx, t \in G\}$, a $\NM$-subgroup of $\M_r(G)$. Thus, the isometric isomorphisms $\vp$ and $\vp^{-1}$ are $\NM$-isomorphisms, so by Theorem \ref{biNM Isomorphism Theorem}, $\vp = \jgt^*$ for some standard isomorphism $(\gamma, \theta)$ of $(F, \omf)$ into $(G, \omg)$.  Applying Lemma \ref{Zadeh Ghar-Zadeh Lemma} (i) as before, for $s \in F$, there is an $\alpha \in \T$ such that  $${1\over \omf(s)}  \gamma(s) \delta_{\theta(s)} = \vp \left({1 \over \omf(s)}  \delta_s\right) = {\alpha\over \omg(\theta(s))}  \delta_{\theta(s)}; \ \ \text{hence} \ \   |\gamma(s)| = \left|{\alpha \omf(s) \over \omg(\theta(s))}\right| = {\omf(s) \over \omg(\theta(s))},$$ as needed. 
\end{proof}

The proofs of the above  theorems for isomorphisms $\vp: \L^1(\omf) \ra \L^1(\omg)$ were independent of Lemma \ref{Zadeh Ghar-Zadeh Lemma} (ii), since  $so_l^\omf - \sigma^\omg$ continuity of the extension of $\vp$ from $\L^1(\omf)$ to $\M(\omf)$ is automatic. We now show that the converse directions of our isomorphism theorems all hold. For bipositive isomorphisms, these are described in \cite{Gha-Zad} as the``easy" directions.   We first establish the  weighted  analogue of \cite[Corollary III.11.20]{Fel-Dor}.

\bp \label{Description of L1(w) in M(w) Prop}  Let  $\mu \in \M(\om)$. Then $\mu \in \L^1(\om)$ if and only if $x \mapsto \mu * \delta_x: G \ra (\M(\om), \|\cdot\|_\om)$ is continuous. 
\ep 

\begin{proof}  Suppose that $x \mapsto \mu * \delta_x$ is continuous. Let $(f_i)$ be a net of positive functions in $C_{00}(G) \subseteq \L^1(\om)$ such that $\s(f_i) \ra {e_G}$ and $\int f_i \, dt = 1$. Let $\ve >0$ and take $i_0$ such that $\| \mu * \delta_t - \mu \|_\om < \ve$ for $t \in \s(f_i)$ and $i \succeq i_0$. Then, for $i \succeq i_0$ and  $\psi \in C_{00}(G)$  with $\| \psi\|_{\infty, \om^{-1}}\leq 1$, 
\beqs |\l \mu *f_i - \mu, \psi \r_\om| & = & \left|   \iint f_i(s^{-1}t) \, d\mu(s) \, \psi(t) \, dt - \int \psi(s) \, d\mu(s)   \right| \\
& = & \left|   \iint f_i(t) \psi(st) \, dt \, d\mu(s) - \iint f_i(t) \psi(s) \, d\mu(s) \, dt  \right| \\
& = & \left|   \int f_i(t) \left( \int \psi(st) \, d\mu(s) - \int  \psi(s) \, d\mu(s)\right) \, dt  \right| \\
& = & \left|   \int f_i(t) \l \mu * \delta_t - \mu, \psi \r_\om  \, dt  \right| \\
& \leq & \int_{\s(f_i)} f_i(t) \| \mu * \delta_t - \mu\|_\om \| \psi\|_{\infty, \om^{-1}}  \, dt \leq \ve,
\eeqs 
where we have used formula (10) in \cite{Sto2} and the Fubini Theorem.  Hence, $\| \mu * f_i - \mu\|_\om \ra 0$; since $\L^1(\om)$ is a closed ideal in $\M(\om)$, $\mu \in \L^1(\om)$. The converse direction is \cite[Lemma 3.1.5]{Zad0}. 
\end{proof} 

In both the weighted and non-weighted cases, the following proposition seems to be new.  
\bp \label{Isomorphisms mapping L1 to L1 Prop} Let $\vp: \M(\omf) \ra \M(\omg)$ be a Banach algebra isomorphism. The following statements are equivalent: 
\bi \item[(i)]  $\vp$ maps $\L^1(\omf)$ isomorphically onto $\L^1(\omg)$;
\item[(ii)]  $\vp$ and $\vp^{-1}$ are $so_l - so_l$ continuous; 
\item[(iii)] $x \mapsto \vp(\delta_x): F \ra (\M(\omg), so_l^\omg)$ and $y \mapsto \vp^{-1}(\delta_y): G \ra (\M(\omf), so_l^\omf)$ are continuous. 
\ei 
\ep

\begin{proof} By Lemmas \ref{Zadeh Ghar-Zadeh Lemma}(ii) and \ref{Extension of isom from L1(w) to M(w) Lemma},  (i) implies (ii) and it is clear that (iii) follows from (ii). For (iii) implies (i), we will use Proposition \ref{Description of L1(w) in M(w) Prop} to observe that continuity of $y \mapsto \vp^{-1}(\delta_y): G \ra (\M(\omf), so_l^\omf)$ implies that $\vp(\L^1(\omf)) \subseteq \L^1(\omg)$: Let $f \in \L^1(\omf)$ and suppose that $y_i \ra y $ in $G$.  Then $\| f*(\vp^{-1} (\delta_{y_i}) - \vp^{-1}(\delta_y))\|_\omf \ra 0$, so $\| \vp(f)*(\delta_{y_i} - \delta_y)\|_\omg \ra 0$; hence $\vp(f) \in \L^1(\omg)$. 
\end{proof}  

\bc The converse directions  of Theorems \ref{biNM Isomorphism Theorem}, \ref{Positive Isom Thm} and \ref{Isometric Isom Thm}  all hold.  
\ec 

\begin{proof} Let $(\gamma, \theta)$ be a standard isomorphism of $(F, \omf)$ onto $(G, \omg)$. The converse statements are satisfied by  the $\NM$-isomorphism $\vp = \jgt^*: \M(\omf) \ra \M(\omg)$ and its inverse $\vp^{-1}= j_{\beta, \theta^{-1}}$ in the various cases by   Proposition \ref{Easy Direction M(w) Prop} (with $\beta$ as defined therein).  If $s_i \ra s$ in $F$, then $\theta(s_i) \ra \theta(s)$ in $G$, so $\delta_{\theta(s_i)} \ra \delta_{\theta(s)}$ $so_l^\omg$; also, $\gamma(s_i) \ra \gamma(s) $ in $\Cx$, so $\vp(\delta_{s_i}) = \gamma(s_i)\delta_{\theta(s_i)} \ra    \gamma(s)\delta_{\theta(s)} = \vp(\delta_s)$ $so_l^\omg$. Hence, $\vp$ and, similarly, $\vp^{-1}$ satisfy condition (iii) of Proposition \ref{Isomorphisms mapping L1 to L1 Prop}, so $\vp$ maps $\L^1(\omf)$ isomorphically onto $\L^1(\omg)$; when $\vp$ is positive or isometric, so is its restriction to $\L^1(\omf)$.
\end{proof} 

That $\jgt^*$ maps $\L^1(\omf)$ onto $\L^1(\omg)$ when $(\gamma, \theta)$ is a positive standard isomorphism is \cite[Theorem 2.2(ii)]{Gha-Zad}, which proved in a different way. 

\br \rm By Corollaries \ref{Sto1 Thm Corollary} and \ref{Positive homs unweighted Corollary} (b), if $\vp: \L^1(F) \ra \M_r(G)$ is a $\NM$-homomorphism, e.g., if $\vp$ is a positive homomorphism, then $\vp$ is automatically contractive. Thus, when $\vp:\L^1(F) \ra \L^1(G)$ or $\vp: \M_r(F) \ra \M_r(G)$ is a $\NM$-isomorphism such that $\vp^{-1}$ is also an $\NM$-isomorphism,  or if $\vp$ is a positive (and therefore bipositive by Theorem \ref{Positive Isom Thm}) isomorphism, then $\vp$ is an isometric isomorphism. By \cite[Example 3.9]{Gha-Zad}, bipositive isomorphisms are not necessarily isomorphic when using  non-trivial weights. 

Corollaries \ref{Sto1 Thm Corollary} and \ref{phi(delta_e) positive Corollary}(b) also imply that if $\vp:\L^1(F) \ra \M_r(G)$ is a $\NM$-homomorphism with $\vp(\delta_{e_F}) = \delta_{e_G}$, then $\vp = \jgt^*$ where $\gamma$ maps $F$ into $\T$. For non-trivial weights, the example of $j_{\gamma^\times, \theta_H}^*$ shows that $\gamma$ does not have to map into $\T$ for such homomorphisms. For an example of a positive isometric  isomorphism $\jgt^*: \M(\omf) \ra \M(\omf)$ such that $\gamma$ does not map into $\T$, take $F = (\R, +)$, $\omf(s) = 2^s$, $\theta: F \ra F: s \mapsto 2s$, $\gamma(s) = \omf(s) / \omf(\theta(s)) = 2^{-s}$. Then $\gamma$ maps onto $(0,\infty)$, but $(\gamma, \theta)$ is a positive linked standard isomorphism of $(F, \omf)$ onto itself, so $\jgt^*: \M(\omf) \ra \M(\omf)$ is a positive isometric isomorphism.   
\er 

\section{$\NM$-homomorphisms $\C F \ra \C G$}

In what follows, $F$ and $G$ are discrete groups.  The purpose of this section is to observe that versions of all the previous results hold for the algebraists' group algebras $\C F$ and $\C G$.  We believe these results are also new. To distinguish it from the pointwise-defined product, we  continue to use the convolution notation $f*g$ for the product of $f$ and $g$ in $\C G$ and, when relevant, we use the $\|\cdot \|_1$-norm on these algebras.  However,  unless explicitly stated otherwise, \it we no longer assume that linear maps on $\C F$ are bounded/continuous with respect to any norm or topology. \rm  We employ the obvious analogues of all definitions previously discussed, dropping all boundedness/continuity conditions (unless explicitly stated otherwise).  In particular: 
\bi \item a subgroup $\Gamma$ of $(\C G, *)$ is a $\NM$-subgroup if $|f*g| = |f|*|g|$  (equivalently, $\|f*g\|_1=\|f\|_1\|g\|_1$) for every $f, g \in \Gamma$;
\item a homomorphism --- not necessarily bounded  ---  $\vp: \C F\ra \C G$ is a $\NM$-homomorphism if  $\Gamma_\vp = \vp(F)$ is a $\NM$-subgroup of $\C G$. 
\ei

 \br \rm Suppose that $\omf$ and $\omg$ are weights on $F$ and $G$.  With respect to the norm $\|f \|_\omg = \|f \omg\|_1$, $\C G $   is a normed algebra, denoted $\C (G, \omg)$, with completion $\ell^1(G, \omg) = \L^1( \omg) = \M(\omg)$ in this, the discrete case. Although we have studied $\NM_{\omg}$-subgroups of $\M(\omg)$ and $\NM_{\omg}$-homomorphisms $\vp: \L^1(\omf) \ra \M(\omg)$, our primary focus has been  on $\NM$-subgroups and $\NM$-homomorphisms. As sets, and algebraically, $\C G = \C (G, \omg)$, so $\vp : \C(F, \omf) \ra \C(G, \omg)$ is a $\NM$-homomorphism if and only if $\vp : \C F \ra \C G $ is a $\NM$-homomorphism (because we are no longer requiring continuity of  $\vp$, which would change things)  if and only if $\Gamma_\vp$ is a $\NM$-subgroup of $\C G$. Thus, in this section in which we are no longer interested in, or required to consider, continuity,  there is little purpose in considering non-trivial weights.  
 \er 
 
 Let $K$ be a finite subgroup  of $G$, $\rho:K \ra \T$ a homomorphism. Then $$m_K = {1\over |K|} \sum_{k\in K} k \qquad \text{and} \qquad \rho m_K = {1\over |K|} \sum_{k\in K} \rho(k) k$$ 
 are (norm-one) idempotents in $\C G$; when $H$ is a subgroup of $G$ and  $K \lhd H$, $m_K \in Z(\C H)$, and when it is also true that $\ker \rho \lhd H$ and $K/\ker \rho \leq Z(H/\ker \rho)$, then $\rho m_K \in Z (\C H)$ (and conversely).

 \bp \label{NM-subgroups in CG}  
 Let  $\Gamma$ be a $\NM$-subgroup of $\C G$ with identity $\iota_\Gamma$. Then: 
\bi \item[(i)] $H = \s(\Gamma) = \bigcup_{f\in \Gamma} \s(f) \leq G$, $\iota_\Gamma = \rho m_K$ for some finite normal subgroup $K$ of $H$ and a homomorphism  $\rho : K \ra \T $ such that \ $\ker \rho \lhd H$ and $K / \ker \rho \leq Z(H/\ker \rho)$;
\item[(ii)] $\Omega_\Gamma = \{(\alpha, t) \in \Cx \times H: \alpha t * \rho m_K \in \Gamma\}$    is a subgroup of   $\Cx \times H$ and
 $\Gamma_{\Cx \times H}= \{ \alpha  t * \rho m_K: (\alpha, t) \in \Cx \times H\}$ is a $\NM$-subgroup of $\C G$ containing $\Gamma$;
\item[(iii)] $ \phi_0: (\alpha, t) \mapsto \alpha t * \rho m_K $
defines a  group homomorphism of $\Omega_\Gamma$ onto $\Gamma$ and of  $\Cx \times H$ onto $\Gamma_{\Cx \times H}$ with kernel  $\Omega_\rho = \{(\rho(k), k): k \in K\}$.
\ei 
If $\Gamma$ is a subgroup  contained in the unit ball of $\C G$,  then $\Gamma$ is a $\NM$-subgroup of $\C G$, and $\C^x$ can  be replaced by $\T $  in the above conditions; if  $\Gamma$ is a positive subgroup of $\C G$, i.e., if each $f \in \Gamma$ is positive, then $\Gamma$ is a $\NM$-subgroup of $\C G$,  $\Cx$ can be replaced by $(0,\infty)$ in the above conditions, and $\iota_\Gamma = m_K$.  \ep  

Though Proposition \ref{NM-subgroups in CG} follows from the results in Section 3, it can be proved directly without considering the topological and measure-theoretic complications inherent in the general case.   

Let $\gamma: F \ra \Cx$, $\theta: F \ra G$ be homomorphisms, $K$ a finite normal subgroup of $H$. The basic homomorphisms act on group elements by 
\medskip 

\begin{tabular}{rl@{\hspace*{15mm}}rl} 
& $\ds M_\gamma^* : \C F \ra \C G: s \mapsto \gamma(s) s,$ & 
& $ \ds  j_\theta^* : \C F \ra \C G: s \mapsto \theta(s), $ \\[2mm]
& $\ds \jgt^*= j_\theta^* \circ M_\gamma^*:  \C F \ra \C G: s \mapsto \gamma(s) \theta(s), $ & 
&   $\ds  S_K^*: \C (H/K) \hookrightarrow \C H: sK \mapsto  s* m_K.$
\end{tabular} 
\medskip 

If we assume that a $\NM$-homomorphism $\vp: \C F \ra \C G$ is continuous, then it has a unique extension to a $\NM$-homomorphism $\vp: \ell^1(F) =\M_r(F) \ra \ell^1(G)= \M_r(G)$. Since the basic homomorphisms map $\C F $ into $\C G$ (etc.), the results in Section 4 immediately yield corresponding results for bounded $\NM$-homomorphisms $\vp: \C F \ra \C G$. However, versions of all statements from Section 4 hold without assuming $\vp: \C F \ra \C G$ is bounded. The proofs from Section 4 can be used as a guide to establishing these results, however the proofs for  homomorphisms $\vp : \C F \ra \C G$ are usually much simpler. We conclude with a sample of  some of these statements. 

\bt   Let $\vp: \C F \ra \C G$ be a mapping. The following statements are equivalent: 
\bi \item[(i)] $\vp$ is a $\NM$-homomorphism; 
\item[(ii)] there is a  subgroup $H$ of $G$, a finite normal subgroup $\Omega_0$ of $\Cx \times H$  and a  homomorphism $\theta:F \ra \Cx \times H / \Omega_0$ such that $\vp = j_{\gamma^\times, \theta_H}^* \circ S_{\Omega_0}^* \circ j_\theta^*$. That is, $\vp$ factors as 
\beq   \label{Main Algebraic Thm Eqn}  \xymatrixrowsep{2pc} \xymatrixcolsep{4pc}
\xymatrix{\C F \ar@<.1ex>[r]^{j_\theta^* \qquad \ \   } & \C(\Cx \times H/\Omega_0)   \ar@{^{(}->}[r]^{S_{\Omega_0}^*} & \C(\Cx \times H)  \ar@<.1ex>[r]^{\quad j_{\gamma^\times, \theta_H}^*}  &   \C G.  } \eeq
\ei
In statement (ii),  we can take $H = \s(\Gamma_\vp)$, $\vp(e_F) = \rho m_K$ and $\Omega_0 = \Omega_\rho$ associated to $\Gamma_\vp$ as in Proposition  \ref{NM-subgroups in CG}.
When $\|\vp(s) \|_1 \leq 1 $ for all $s \in F$ --- in which case $\vp$ is contractive and automatically a $\NM$-homomorphism --- $\Cx$ can be replaced by $\T$ in (\ref{Main Algebraic Thm Eqn}).  
\et 

\bc \label{Algebraic Cor}1.  Let $\vp: \C F \ra \C G$ be a $\NM$-homomorphism with $H = \s(\Gamma_\vp)$ and $\iota_{\Gamma_\vp}= \rho m_K$, (as in Proposition \ref{NM-subgroups in CG}). If $\rho$ extends to a homomorphism $\rho_H: H \ra \T$, then there are homomorphisms $\gamma: F \ra \Cx$ and $\theta : F \ra H/K$ such that $\vp= M_{\rho_H}^* \circ S_K^* \circ j_{\gamma, \theta}^*$, i.e., $\vp $ factors as
\beq    \xymatrixrowsep{2pc} \xymatrixcolsep{4pc}
\xymatrix{\C F \ar@<.1ex>[r]^{j_{\gamma,\theta}^* \quad \ \   } & \C(H/K)   \ar@{^{(}->}[r]^{S_K^*} & \C H  \ar@<.1ex>[r]^{ M_{\rho_H}^*\qquad }  &   \C H \subseteq \C G.  } \eeq
2.  A mapping $\vp: \C F \ra \C G$ is a $\NM$-homomorphism with $\vp(e_F) \geq 0 $ if and only if there  is a  subgroup $H$ of $G$, a finite normal subgroup $K$ of $H$, and homomorphisms $\gamma: F \ra \Cx$ and $\theta : F \ra H/K$ such that $\vp= S_K^* \circ j_{\gamma, \theta}^*$.

\smallskip 

\noindent  3.   A mapping $\vp: \C F \ra \C G$ is a positive homomorphism if and only if there  is a  subgroup $H$ of $G$, a finite normal subgroup $K$ of $H$, and homomorphisms $\gamma: F \ra (0, \infty)$ and $\theta : F \ra H/K$ such that $\vp= S_K^* \circ j_{\gamma, \theta}^*$.
\ec 

Algebraic versions of the isomorphism theorems from Section 4.3 can be deduced from Corollary \ref{Algebraic Cor}, however these are easy to establish without any additional machinery.

\noindent{\sc Department of Combinatorics and Optimization, University of Waterloo, Waterloo, ON, N2L 3G1, Canada}

\noindent email address: {\tt mekroeker@uwaterloo.ca }

\bigskip

\noindent {\sc Department of Mathematics, University of Manitoba, Winnipeg, Manitoba, Canada R3T 2N2} 

\noindent email address: {\tt  stephe47@myumanitoba.ca}

\bigskip 

\noindent {\sc Department of Mathematics and Statistics, University
of Winnipeg, 515 Portage Avenue, Winnipeg, MB, R3B 2E9, Canada }

\noindent email address: {\tt r.stokke@uwinnipeg.ca}

\bigskip 

\noindent{\sc Department of Computer Science, University of Calgary, 2500 University Drive NW 
Calgary, Alberta, Canada, T2N 1N4} 

\noindent email address: {\tt  randy.yee1@ucalgary.ca}

\end{document}

%% file: Macros.tex
\usepackage{amssymb}

\newtheorem{thm}{Theorem}[section]

\newtheorem{dfn}[thm]{Definition}

\newtheorem{lem}[thm]{Lemma}

\newtheorem{prop}[thm]{Proposition}

\newtheorem{remark}[thm]{Remark}

\newtheorem{cor}[thm]{Corollary}

\newtheorem{ex}[thm]{Example}

\newtheorem{question}[thm]{Question}

\def\ve{\varepsilon}
\def\sq{{\scriptscriptstyle \square}}

\def\bq{\begin{question}}
\def\bt{\begin{thm}}
\def\bp{\begin{prop}}
\def\blem{\begin{lem}}
\def\bd{\begin{dfn}}
\def\br{\begin{remark}}
\def\bc{\begin{cor}}
\def\bex{\begin{ex}}
\def\beqs{\begin{eqnarray*}}
\def\beq{\begin{eqnarray}}
\def\bi{\begin{itemize}}

\def\eq{\end{question}}
\def\et{\end{thm}}
\def\ep{\end{prop}}
\def\elem{\end{lem}}
\def\ed{\end{dfn}}
\def\er{\end{remark}}
\def\ec{\end{cor}}
\def\eex{\end{ex}}
\def\eeqs{\end{eqnarray*}}
\def\eeq{\end{eqnarray}}
\def\ei{\end{itemize}}

\def\ds{\displaystyle}

\def\c{\cdot}

\def\ov{\overline}
\def\r{\rangle}
\def\l{\langle}

\def\w*{w^*-w^*}

\def\ra{\rightarrow}

\def\C{\cdot}

\def\vp{\varphi}
 
\def\bs{\backslash}

\def\om{\omega} 
\def\M{{\cal M}}
 
\def\fS{\mathfrak{S}}

\def\L{{\cal L}}

\def\NMO{\rm{NM}_\om} 
\def\NM{\rm{NM}} 
\def\ig{\iota_\Gamma}

\makeatletter
\def\moverlay{\mathpalette\mov@rlay}
\def\mov@rlay#1#2{\leavevmode\vtop{%
   \baselineskip\z@skip \lineskiplimit-\maxdimen
   \ialign{\hfil$\m@th#1##$\hfil\cr#2\crcr}}}
\newcommand{\charfusion}[3][\mathord]{
    #1{\ifx#1\mathop\vphantom{#2}\fi
        \mathpalette\mov@rlay{#2\cr#3}
      }
    \ifx#1\mathop\expandafter\displaylimits\fi}
\makeatother

%% file: NM_homomorphisms_of_Beurling_Algebras_3.bbl
\begin{thebibliography}{99}


\bibitem{Coh} P. J. Cohen, On homomorphisms of group algebras, {\it
Amer. J. Math} 82 (1960), 213-226.


\bibitem{Dal-Lau} H.G. Dales and A. T.-M. Lau, The second duals of Beurling algebras, {\it  Mem. Amer. Math. Soc.}  177 (2005), no. 836.


\bibitem{Fel-Dor}  J.M.G.  Fell and R.S. Doran, {\it Representations of $*$-algebras,
locally compact groups, and Banach $*$-algebraic bundles  Vol.
1},  Pure and Applied Mathematics, 125. Academic Press, Inc., 1988.

\bibitem{Fil-Sal} M.  Filali and P.  Salmi, Topological centres of weighted convolution algebras, \it  J. Funct. Anal. \rm  278 (2020), no. 11, 22 pages.

\bibitem{Gha1}  F. Ghahramani, Weighted group algebra as  an ideal in its second dual space, \it Proc. Amer. Math.  Soc. \rm  90 (1984), no. 1, 71-76.


\bibitem{Gha2} F. Ghahramani, Compact elements of weighted group algebras, {\it Pacific J. Math.} 113 (1984), no. 1, 77-84.

\bibitem{Gha-Lau1} F.  Ghahramani and A.T.  Lau,  Isometric isomorphisms between the
second conjugate algebras of group algebras, {\it  Bull. London
Math. Soc.}  20 (1988),  no. 4, 342--344.

\bibitem{Gha-Lau-Lo} F. Ghahramani, A.T. Lau, and V. Losert,  Isometric isomorphisms between Banach algebras related to locally compact groups,
{\it Trans. Amer. Math. Soc.} 321 (1990), no. 1, 273-283.

\bibitem{Gha-Lau2} F. Ghahramani and A.T-M. Lau, Isomorphisms and
multipliers on second dual algebras of Banach algebras, {\it   Math.
Proc. Cambridge Philos. Soc.}  111  (1992),  no. 1, 161-168.

\bibitem{Gha-Zad} F.  Ghahramani and S.  Zadeh, Bipositive isomorphisms between Beurling algebras and between their second dual algebras,  {\it Canad. J. Math.}  69 (2017), no. 1, 3-20


\bibitem{Gli}  I. Glicksberg,  Homomorphisms of certain algebras of measures, {\it Pacific J. Math.} 10 (1960), 167-191.

\bibitem{Gre} F.P. Greenleaf,  Norm decreasing homomorphisms of
 group algebras,  {\it Pacific J. Math.} 15 (1965), 1187-1219.
 
 

\bibitem{Joh} B.E. Johnson, Isometric isomorphisms of measure algebras, {\it Proc. Amer. Math. Soc.} 15 (1964), 186-188.

\bibitem{Kan} E. Kaniuth,  {\it A course in commutative Banach algebras}, Springer, New York, 2009.

\bibitem{Kal-Woo} N.J. Kalton and G.V. Wood,  Homomorphisms of group algebras with
norm less than $\sqrt{2}$,  {\it Pacific J. Math.}   62  (1976), no.
2, 439--460.


\bibitem{Kaw} Y. Kawada, On the group ring of a topological group, {\it Math. Japonicae} 1(1948), 1-5



\bibitem{Lau-McK} A. T.-M.  Lau and K. McKennon, Isomorphisms of locally
compact groups and Banach algebras,   {\it Proc. Amer. Math. Soc. }
79 (1980), no. 1, 55--58.

\bibitem{Neu} M. Neufang, On the topological centre problem for weighted convolution algebras and semigroup compactifications, {\it Proc. Amer. Math. Soc.} 136 (2008), no. 5, 1831-1839. 

\bibitem{Num} K. Numakura,  On bicompact semigroups, {\it Math. J. Okayama Univ.} 1 (1952), 99-108. 



\bibitem{Pym} J. S. Pym,  Idempotent measures on semigroups, {\em Pacific J. Math.} 12 (1962), 685-698.


\bibitem{Rei-Ste}  H. Reiter and J.D. Stegeman, {\it
Classical harmonic analysis and locally compact groups}, 
second ed., London Math. Soc. Monographs,  Volume 22, Clarendon Press,  Oxford,  2000.



\bibitem{Sam} E.  Samei,  Weak amenability and 2-weak amenability of Beurling algebras, {\it J. Math. Anal. Appl.} 346 (2008), no. 2, 451-467.

\bibitem{She-Zha} V. Shepelska and Y.  Zhang, Non-weakly amenable Beurling algebras, {\it Indiana Univ. Math. J.} 67 (2018), no. 1, 119-150.


\bibitem{Spr} N. Spronk, Commuting contractive idempotents in measure algebras, {\it Ann. Funct. Anal.} 7 (2016), no. 1, 136-149.
 
   \bibitem{Sto1} R. Stokke, Homomorphisms of convolution algebras, {\it J. Funct. Anal.} 261 (2011), no. 12, 3665-3695.
   
   
      \bibitem{Sto2} R. Stokke, On Beurling measure algebras, {\it  Comment. Math. Univ. Carolin.}, in press.
      
 \bibitem{Str} R.S. Strichartz,  Isometric isomorphisms of measure
algebras, {\it Pacific J. Math.},  15  (1965), 315-317.

\bibitem{Wal} M. Walter,  $W\sp{*} $-algebras and nonabelian harmonic
analysis, {\it  J. Funct. Anal.} 11 (1972), 17-38.

 \bibitem{Wen1} J.G.Wendel, On isometric isomorphisms of group algebras, {\it Pacific J. Math.} 1(1951), 305-311.
 
 \bibitem{Wen2} J.G. Wendel,  Left centralizers and isomorphisms of group algebras. {\it Pacific J. Math.} 2 (1952), 251–261.
 
\bibitem{Zad0} S. Zadeh, Isomorphisms of Banach algebras associated with locally compact groups, PhD thesis, University of Manitoba, Canada, 2015. 

\bibitem{Zad} S. Zadeh,  Isometric isomorphisms of Beurling algebras, {\it  J. Math. Anal. Appl.} 438 (2016), 
 1-13. 
 
 
 

\end{thebibliography}
